\documentclass[leqno]{amsart}

\usepackage[english]{babel} 									
\usepackage[T1]{fontenc}							
\usepackage[utf8]{inputenx}
				
\usepackage{mathtools}								
\usepackage{empheq}									
\usepackage{enumitem}								

\usepackage{amssymb}								
\usepackage{mathrsfs}								
\usepackage{bm}                                     
				
\usepackage{verbatim}                               
\usepackage{iflang}									
\usepackage{xifthen}								

\usepackage[final]{pdfpages}						
\usepackage{todonotes}								
\usepackage{aliascnt}								
\usepackage{needspace}								
\usepackage{refcount}

\usepackage{subcaption}			
\usepackage{placeins} 			
\usepackage{grffile}			
\usepackage[all]{xy}			

\usepackage{transparent}
\usepackage{color}				

\usepackage{tikz}               
\usetikzlibrary{math,arrows.meta, cd}
\DeclareGraphicsExtensions{.pdf,.png,.jpg}

\usepackage{graphicx}			
 
\usepackage{diagbox} 

\usepackage[babel]{csquotes}							

\usepackage{nameref}									
\usepackage[colorlinks	=	true,						
          	 	linkcolor	=	red,
            	urlcolor	=	red,
            	citecolor	=	red]%
            	{hyperref}
\usepackage{bookmark}									
\usepackage{cleveref}									
\numberwithin{equation}{section}

\newtheorem{theorem}{Theorem}[section]

\newtheorem{corollary}[theorem]{Corollary}
\newtheorem{lemma}[theorem]{Lemma}
\newtheorem{proposition}[theorem]{Proposition}

\newtheorem*{theorem*}{Theorem}

{\bf}{\it}

\newtheorem{definition}[theorem]{Definition}

\theoremstyle{remark}
\newtheorem{remark}[theorem]{Remark}

\newtheorem*{remark*}{Remark}

\theoremstyle{remark}

\newcommand{\N}{\mathbb{N}}
\renewcommand{\S}{\mathbb{S}}

\newcommand{\R}{\mathbb{R}}

\allowdisplaybreaks

\newcommand{\apmd}[2][]{							
	\ifthenelse{\equal{#1}{}}%
					{ \operatorname{N}_{#2}	}%
					{ \operatorname{N}_{#1,#2} 	}}

\newcommand{\aint}[2][]{
	\ifthenelse{\equal{#1}{}}%
					{%
\mathchoice%
      {\mathop{\kern 0.2em\vrule width 0.6em height 0.69678ex depth -0.58065ex
              \kern -0.8em \intop}\nolimits_{\kern -0.45em#2}^{#1}}%
      {\mathop{\kern 0.1em\vrule width 0.5em height 0.69678ex depth -0.60387ex
              \kern -0.6em \intop}\nolimits_{#2}^{#1}}%
      {\mathop{\kern 0.1em\vrule width 0.5em height 0.69678ex depth -0.60387ex
              \kern -0.6em \intop}\nolimits_{#2}^{#1}}%
      {\mathop{\kern 0.1em\vrule width 0.5em height 0.69678ex depth -0.60387ex
              \kern -0.6em \intop}\nolimits_{#2}^{#1}}}%
					{%
\mathchoice%
      {\mathop{\kern 0.2em\vrule width 0.6em height 0.69678ex depth -0.58065ex                                              
              \kern -0.8em \intop}\nolimits_{\kern -0.45em#1}^{#2}}%
      {\mathop{\kern 0.1em\vrule width 0.5em height 0.69678ex depth -0.60387ex
              \kern -0.6em \intop}\nolimits_{#1}^{#2}}%
      {\mathop{\kern 0.1em\vrule width 0.5em height 0.69678ex depth -0.60387ex
              \kern -0.6em \intop}\nolimits_{#1}^{#2}}%
      {\mathop{\kern 0.1em\vrule width 0.5em height 0.69678ex depth -0.60387ex
              \kern -0.6em \intop}\nolimits_{#1}^{#2}}}}

\newcommand{\vol}{\mathrm{vol}}

\usepackage{thmtools}
\usepackage{thm-restate}

\newcommand{\C}{\mathbb{C}}

\newcommand{\SL}{\mathrm{SL}}

\newcommand{\assoc}{\mathrm{assoc}}
\newcommand{\coassoc}{\mathrm{coassoc}}
\newcommand{\Cayley}{\mathrm{Cayley}}

\begin{document}
\title{Singular sets of conformal curves are discrete}
\author{Toni Ikonen}

\address{Department of Mathematics and Statistics\\ P. O. Box 68 (Pietari Kalmin katu 5) \\ FI-00014 University of Helsinki\\ Finland}
\email{toni.ikonen@helsinki.fi}

\thanks{The author was supported by the Swiss National Science Foundation grant 212867 and the Research Council of Finland, project number 374557. }

\begin{abstract}
   We prove that the differential of a non-constant conformal curve between Euclidean spaces has full rank outside a discrete set. In particular, conformal curves are real-analytic outside a discrete set. We provide examples of entire associative and coassociative conformal curves with a non-empty singular set. We show that every non-constant Special Lagrangian conformal curve is a real-analytic conformal immersion. We also prove that a non-constant conformal curve has discrete fibers.
\end{abstract}

\maketitle

\section{Introduction}
\subsection{Main results}
It is well-known that a holomorphic curve has a full-rank differential outside a discrete set. Similarly, Liouville's theorem (1850) states that a conformal mapping from a Euclidean domain $\Omega \subset \R^n$, for $n \geq 3$, into $\R^n$ is the restriction of a Möbius transformation and thus the differential has full rank everywhere. These theorems extend to conformal curves as follows.
\begin{theorem}\label{theorem-full-rank-set}
    The singular set of a non-constant conformal curve $\R^n \supset \Omega \to \R^m$ is discrete for every $3 \leq n \leq m$.
\end{theorem}
In the equidimensional case $n = m$, our proof in fact shows that the singular set is empty, without appealing to Liouville's theorem.

To define conformal curves, we recall the definition of a constant-coefficient calibration $\omega \in \Lambda^{n}\mathbb{R}^m$ for some $2 \leq n \leq m$. Here $\Lambda^n \R^m$ can be canonically identified with differential $n$-forms that have constant coefficients in the standard basis in which case $\omega$ being a calibration simply refers to the comass normalization
\begin{equation*}
    \| \omega \|
    =
    \max
    \left\{
        \omega( v_1 \wedge \dots \wedge v_n )
        \colon
        ( v_1, \dots, v_n )
        \text{ orthonormal}
    \right\}
    =
    1.
\end{equation*}
Calibrations were introduced by Harvey and Lawson in their foundational paper~\cite{harvey-lawson-1982-calibrated-geometries} as a source of new minimal surfaces.

A \emph{conformal curve} $F \colon \R^n \supset \Omega \to ( \R^m, \omega )$ is a map that satisfies
\begin{equation}\label{eq:quasiregular:curve}
    F \in C^{1}( \Omega, \mathbb{R}^m )
    \quad\text{and}\quad
    \| DF \|^n = \star F^{*}\omega
    \quad\text{everywhere},
\end{equation}
where $\|DF\|$ is the operator norm of the differential, $F^{*}\omega$ is the standard pullback, and $\star$ the Hodge star operator. Here \eqref{eq:quasiregular:curve} can be equivalently expressed as
\begin{equation*}
    F \in C^{1}( \Omega, \mathbb{R}^m )
    \quad\text{and}\quad
    \frac{ 1 }{ n^{ \frac{n}{2} } } | DF |^n = \star F^{*}\omega
    \quad\text{everywhere},
\end{equation*}
for the Hilbert--Schmidt norm $|DF|$. For an angle-preserving embedding $F \in C^1( \Omega, \R^m )$, the equality \eqref{eq:quasiregular:curve} holds if and only if the submanifold $F( \Omega ) \subset \R^m$ is $\omega$-calibrated. In particular, $F( \Omega )$ is a minimal submanifold.

When $\omega \in \Lambda^{2} \C^{k}$ is a symplectic form, the associated conformal curves are exactly holomorphic curves. Similarly, if $\omega \in \Lambda^n \R^n$ is a calibration, then $\omega$ is either the standard volume form $\vol_{\R^n}$ or $-\vol_{\R^n}$, and the associated conformal curves are precisely (orientation-preserving or -reversing, respectively) conformal maps $\R^n \supset \Omega \to \R^n$.  The above generalization of these maps was introduced independently by Smith \cite{smith-2011-a-theory-of-multiholomorphic-maps} and by Pankka \cite{pankka-2020-quasiregular-curves}. The research on conformal curves has been active in recent years. Indeed, conformal curves have also been studied under the name \emph{Smith maps} in \cite{smith-2011-a-theory-of-multiholomorphic-maps,cheng-karigiannis-madnick-2020-bubble-tree-convergence-of-conformally-cross-product-preserving-maps,iliashenko-karigiannis-2025-a-special-class-of-k-harmonic-maps-inducing-calibrated-fibrations,broder-iliashenko-madnick-2025-hyperbolicity-and-schwarz-lemmas-in-calibrated-geometry,iliashenko-karigiannis-madnick-2026-montels-theorem}, under the name \emph{calibrated curve} in \cite{heikkila-pankka-prywes-2023-quasiregular-curves-of-small-distortion-in-product-manifolds}, and as conformal curves in \cite{ikonen-pankka-2024-liouvilles-theorem-in-calibrated-geometries,ikonen-2025-entire-conformal-curves,ikonen-pim-2026-a-rescaling-principle-for-quasiregular-curves-with-applications-to-hyperbolicity,iliashenko-2026-branch-obstructions-for-area-minimizing-currents-and-flat-critical-points-of-weakly-conformal-maps}.

We recall that a conformal curve is $n$-harmonic, cf. \cite{cheng-karigiannis-madnick-2020-bubble-tree-convergence-of-conformally-cross-product-preserving-maps,heikkila-pankka-prywes-2023-quasiregular-curves-of-small-distortion-in-product-manifolds}, and thus $C^{1,\alpha}_{loc}$-regular for some $\alpha \in (0,1)$ \cite{uhlenbeck-1977-regularity-for-a-class-of-non-linear-elliptic-systems,hardt-lin-1987-mappings-minimizing-the-lp-norm-of-the-gradient}. By standard elliptic regularity, the curves are real-analytic in the \emph{full-rank set} $U_F = \{ |DF|> 0 \}$. The structure of the \emph{singular set} $S_F = \{ |DF| = 0 \}$ was previously not well-understood. In particular, Theorem~\ref{theorem-full-rank-set} answers open questions by Heikkilä--Pankka--Prywes, Iliashenko--Karigiannis, and Ikonen--Pankka \cite{heikkila-pankka-prywes-2023-quasiregular-curves-of-small-distortion-in-product-manifolds,iliashenko-karigiannis-2025-a-special-class-of-k-harmonic-maps-inducing-calibrated-fibrations,ikonen-pankka-2024-liouvilles-theorem-in-calibrated-geometries} regarding the structure of the singular set $S_F$.

We next present further results. As an immediate consequence of Theorem~\ref{theorem-full-rank-set}, conformal curves satisfy a unique continuation principle. We state this separately since its direct proof is an ingredient in the proof of Theorem~\ref{theorem-full-rank-set}.
\begin{theorem}\label{theorem-unique-continuation-principle-open-set}
    A conformal curve $\R^n \supset \Omega \to \R^m$ that is constant on a subdomain is constant.
\end{theorem}
We also prove that the fibers of conformal curves have a rigid structure.
\begin{theorem}\label{theorem-discreteness}
    Each fiber of a non-constant conformal curve $\R^n \supset \Omega \to \R^m$ is discrete.
\end{theorem}
The proofs of these theorems also reveal a strong connection to singularities of calibrated cones. In fact, Theorems~\ref{theorem-full-rank-set} and \ref{theorem-discreteness} are sharp due to the following examples.
\begin{theorem}\label{theorem-singular-set-can-be-non-empty}
    There exist radially homogeneous conformal curves
    $$\R^3 \to ( \R^7, \omega_{\assoc} ) \in C^2( \R^3, \R^7 ) \setminus C^{3}( \R^3, \R^7 )$$
    and
    $$\R^4 \to ( \R^7, \omega_{\coassoc} ) \in C^{\omega}( \R^4, \R^7 )$$ with singular set $\{0\}$. The associative curve is a topological embedding and, outside the origin, the cardinality of the fibers for the coassociative curve is 24.
\end{theorem}
The associative conformal curve $\R^3 \to ( \R^7, \omega_{\assoc} )$ shows that conformal curves that are topological immersions (or even embeddings) can have singularities. This contrasts the two-dimensional theory of holomorphic curves, Liouville's theorem, and Ferrand's regularity theorem for conformal mappings on smooth manifolds \cite{lelong-ferrand-1976-geometrical-interpretations-of-scalar-curvature-and-regularity-of-conformal-homeomorphisms}. The coassociative example illustrates that Theorem~\ref{theorem-discreteness} is not an immediate consequence of Theorem~\ref{theorem-full-rank-set}.

Despite the examples from Theorem~\ref{theorem-singular-set-can-be-non-empty}, the following theorem shows that many conformal curves are indeed real-analytic conformal immersions.
\begin{theorem}\label{theorem-special-lagrangian-spheres}
    For $n \geq 3$, every non-constant conformal curve $\R^n \supset \Omega \to ( \R^{2n}, \omega_{\SL} )$ for the Special Lagrangian calibration $\omega_{\SL}$ is a real-analytic conformal immersion.
\end{theorem}
We recall that there are Special Lagrangian cones that are not real-analytic, cf. \cite[page 97]{harvey-lawson-1982-calibrated-geometries}. In particular, due to Theorem~\ref{theorem-special-lagrangian-spheres}, the singularities of such cones cannot be modelled by Special Lagrangian conformal curves. We note that examples of Special Lagrangian curves were described in \cite{ikonen-pankka-2024-liouvilles-theorem-in-calibrated-geometries}. 

Combining Theorem~\ref{theorem-special-lagrangian-spheres} with the rigidity theorem \cite[Theorem 1.7]{ikonen-pankka-2024-liouvilles-theorem-in-calibrated-geometries} leads to the following corollary.
\begin{corollary}\label{corollary-low-dimensional-curves-are-real-analytic}
    Every non-constant conformal curve $\R^{n} \supset \Omega \to \R^m$ is a real-analytic immersion for every $3 \leq n \leq m \leq 6$.
\end{corollary}
In particular, the examples from Theorem~\ref{theorem-singular-set-can-be-non-empty} are dimensionally sharp.

Since Theorems~\ref{theorem-full-rank-set}--\ref{theorem-discreteness} are well-known when $n = 2$, we solely focus on the higher-dimensional case, $n \geq 3$, in this manuscript.

\subsection{Outline of the manuscript}
Section~\ref{section-preliminaries} contains our preliminaries. More precisely, in Section~\ref{section-basic-notation}, we introduce our basic notations. In Section~\ref{section-natural-differential}, we recall a standard proposition on the first-order regularity of the natural differential $V_F = |DF|^{ \frac{n-2}{2} } DF$ of $n$-harmonic maps: $V_F \in W^{1,2(1+\delta)}_{loc}(\Omega)$ for some $\delta > 0$.

Section~\ref{section-higher-order-regularity} starts with a key second-order Sobolev regularity result for $|V_F|^2$ when $F$ is a conformal curve. That is, we prove that $|V_F|^2 \in W^{2,1+\delta}_{loc}( \Omega )$, see Theorem~\ref{theorem-higher-order-regularity}, with the parameter $\delta$ from Section~\ref{section-natural-differential}. This regularity leads to a conservation law for conformal curves. Indeed, the matrix field
    \begin{align}\label{equation-weitzenbock-bochner-formula-matrix-version}
        \mathcal{T}^F 
        = 
        \frac{n-2}{n^2}\left( D^{2}( |V_F|^2 ) - \Delta( |V_F|^2 ) I_{ n \times n } \right)
    \end{align}
is locally integrable and it is distributionally divergence-free. The analysis of $\mathcal{T}^F$ leads to a key result in Section~\ref{section-frequency-function-for-conformal-curves}. Indeed, we define the \emph{frequency function} 
\begin{align*}
    \mathcal{N}_{F}(x,r)
    =
    \frac{r}{n}
    \frac{ \int_{ B(x,r) } \Delta( |V_F|^2 ) \,dz }{ \int_{ \partial B(x,r) } |V_F|^2 \,ds }
\end{align*}
whenever the denominator is positive and zero otherwise. Then we prove the following theorem.
\begin{theorem}\label{theorem-frequency-function-monotonicity}
    Let $3 \leq n \leq m$ and let $F \colon \R^n \supset \Omega \to ( \R^m, \omega )$ be a conformal curve. Whenever $x \in \Omega$, the function $r \mapsto \mathcal{N}_{F}(x,r)$ is non-decreasing for each $0<r<\mathrm{dist}(x, \R^n \setminus \Omega )$.
\end{theorem}
In particular, the subsequent definition is well-posed.
\begin{definition}[The pointwise frequency]
    Let $3 \leq n \leq m$ and let $F \colon \R^n \supset \Omega \to ( \R^m, \omega )$ be a conformal curve. The \emph{pointwise frequency} of $F$ at $x$ is
    \begin{align}\label{equation-introduction-pointwise-frequency}
        \mathcal{N}_{F}(x) = \lim_{ r \to 0^{+} } \mathcal{N}_{F}(x,r).
    \end{align}
\end{definition}  
The formula for $\mathcal{N}_{F}(x,r)$ is inspired by the Poisson--Jensen formula for subharmonic functions, Almgren's frequency formula for harmonic functions \cite{almgren-1979-dirichlets-problem-for-multiple-valued-functions-and-the-regularity-of-mass-minimizing-integral-currents}, and the Weitzenböck--Bochner type formula \cite{takeuchi-1994-some-conformal-properties-of-p-harmonic-maps-and-a-regularity-for-sphere-valued-p-harmonic-maps}.

In Section~\ref{section-radially-homogeneous-conformal-curves}, we analyze conformal curves for which $r \mapsto \mathcal{N}_{F}(x,r)$ is a constant $\mathcal{N}_0 \in [0,\infty)$ using \eqref{equation-introduction-pointwise-frequency} and the monotonicity result, Theorem~\ref{theorem-frequency-function-monotonicity}. Such curves are described by the following theorem.
\begin{theorem}\label{theorem-introduction-the-main-frequency-structure-result}
    Let $3 \leq n \leq m$ and let $F \colon \R^n \supset \Omega \to ( \R^m, \omega )$ be a conformal curve. Let $\mathcal{N}_0 \in [0,\infty)$. If $\mathcal{N}_{F}(x,r) = \mathcal{N}_0$ for every $0 < r < \mathrm{dist}(x, \R^n \setminus \Omega )$, then there exists a real-analytic $h \colon \S^{n-1} \to \R^{m}$ such that
    \begin{align}\label{equation-introduction-frequency-constant-functions}
        F(x+rv) = F(x) + r^{ \mathcal{N}_0 + 1 } h(v)
        \quad\text{whenever $x+rv \in \Omega$.}
    \end{align}
    Furthermore, there exists a constant $C \geq 0$ such that the following holds.
    \begin{enumerate}
        \item $h$ maps into a sphere of radius $C$;
        \item the pullback metric satisfies $h^{*}g_{euc} = C^2 \left( \mathcal{N}_{0}+1\right)^2 g_{\S^{n-1}}$;
        \item the spherical Laplacian satisfies $\Delta_{ \S^{n-1} } h = -(n-1) \left( \mathcal{N}_{0} + 1 \right)^2 h$.
    \end{enumerate}
\end{theorem}
    Such $h$ were analyzed by Takahashi \cite{takahashi-1966-minimal-immersions-of-riemannian-manifolds} and do Carmo--Wallach \cite{do-carmo-wallach-1971-minimal-immersions-of-spheres-into-spheres} in the context of minimal isometric immersions of spheres into spheres.

    We also prove a blow-up theorem in Section~\ref{section-radially-homogeneous-conformal-curves}. More precisely, consider the normalizing constant
    \begin{align}\label{equation-correct-normalization-used-for-blow-ups}
        c_{r}(x)
        &=
        r
        \left( \aint{ \partial B(x,r) } \frac{ |DF|^{n} }{ n^{ \frac{n}{2} } } \,ds \right)^{ \frac{1}{n} }.
    \end{align}
    A subtle consequence of Theorem~\ref{theorem-unique-continuation-principle-open-set} is that $c_r(x) > 0$ (formulated in Lemma~\ref{lemma-consequence-of-monotonicity}). This allows us to consider the following normalization
    \begin{align}\label{equation-correct-normalization-blow-up}
        F_r(y) = \frac{ F(x+ry) - F(x) }{ c_r(x) }
        \quad\text{and}\quad\Omega_{r} \coloneqq \{ y \in \R^n \colon x+ry \in \Omega \}.
    \end{align}
    The subsequent theorem describes blow-ups at each $x \in \Omega$.
    \begin{theorem}[Radially homogeneous blow-ups]\label{theorem-normalized-blow-up}
        Let $3 \leq n \leq m$ and let $F \colon \R^n \supset \Omega \to ( \R^m, \omega )$ be a non-constant conformal curve. If $(r_j)$ is a sequence of radii that satisfies $r_j \rightarrow 0^{+}$, then, up to passing to a subsequence, $$( F_{r_j} \colon \Omega_{r_j} \rightarrow (\R^m,\omega) )$$ converges uniformly and in $C^1$ on compact sets of $\R^n$ to a non-constant conformal curve $F_{\infty} \colon \R^n \to ( \R^m, \omega )$ such that
        \begin{align*}
            \mathcal{N}_{F_\infty}( 0, r ) = \mathcal{N}_{F}(x)
            \quad\text{and}\quad
            F_\infty(r v) = r^{ \mathcal{N}_F(x) + 1 } h_\infty(v),
            \quad r>0, \quad v \in \S^{n-1},
        \end{align*}
        for an isometric immersion $h_\infty \colon \S^{n-1} \to \R^m$ satisfying
        \begin{enumerate}
            \item $|h_\infty| = \frac{ 1 }{ \mathcal{N}_{F}(x) + 1 }$;
            \item $\Delta_{ \S^{n-1} } h_\infty = -(n-1) \left( \mathcal{N}_{F}(x) + 1 \right)^2 h_\infty$.
        \end{enumerate}
    \end{theorem}
    A key point of Theorem~\ref{theorem-normalized-blow-up} is that the subsequential limits are described using Theorem~\ref{theorem-introduction-the-main-frequency-structure-result}. The normalizing constant \eqref{equation-correct-normalization-used-for-blow-ups} is chosen to guarantee that $h_{\infty}$ is an isometric immersion.
    
    In Section~\ref{section-proof-of-the-main-results}, we prove the main results. For instance, Theorem~\ref{theorem-full-rank-set} is immediate from the blow-up theorem, Theorem~\ref{theorem-normalized-blow-up}: an accumulation point of singular points of a non-constant conformal curve produces a singular point on the unit sphere for the radially homogeneous blow-up, contradicting Theorem~\ref{theorem-the-main-frequency-structure-result}. A similar blow-up argument underlies Theorem~\ref{theorem-discreteness}. Theorem~\ref{theorem-normalized-blow-up} leads to a direct relationship between the singular set of non-constant conformal curves and radially homogeneous conformal curves.
\begin{theorem}\label{theorem-singular-set-is-empty-if-and-only-if-radial-homogeneity}
    Let $3 \leq n \leq m$ and let $\omega \in \Lambda^n \R^m$ be a calibration. Then the following are equivalent.
    \begin{enumerate}
        \item The singular set of every non-constant conformal curve $\R^n \to ( \R^m, \omega )$ is empty.
        \item Every radially homogeneous conformal curve $\R^n \to ( \R^m, \omega )$ is linear.
        \item Every $\omega$-calibrated cone over an immersed isometric sphere $\S^{n-1}(0,c) \rightarrow \S^{m-1}$ is an $n$-plane.
    \end{enumerate}
\end{theorem}
To clarify the statement, the cone over a minimal isometric immersion $$h \colon \S^{n-1}(0,c) \rightarrow \S^{m-1}$$ is defined by the spherical coordinate representation $rv \mapsto r h(cv)$, or, equivalently, by the conformal map $rv \mapsto r^{c} h(cv)$.

\section{Preliminaries}\label{section-preliminaries}

\subsection{Standard notation}\label{section-basic-notation}
Our notation is fairly standard. Open Euclidean balls are denoted by $B(x,r)$, closed ones by $\overline{B}(x,r)$. The spheres in $\R^n$ are denoted by $\S^{n-1}(x,r)$. For the unit ball and sphere, we use the notations $B$, $\overline{B}$, and $\S$. In this manuscript, $\Omega \subset \R^n$ denotes a non-empty domain, i.e., a non-empty connected and open set. The relative closure of $A \subset \Omega$ is the set $\Omega \cap \overline{A}$, where $\overline{A}$ denotes the closure of $A$ in $\R^n$. The interior of a set $A \subset \R^n$ is denoted by $\mathrm{int}(A)$. The boundary is denoted by $\partial A$. In some of our estimates, the distance from $x \in \Omega$ to the complement of $\Omega$ plays a role. The distance is denoted by $r(x) = \mathrm{dist}(x, \R^n \setminus \Omega)$.

The Lebesgue measure of $B$ is $\omega_n$ and the unit sphere $\S$ has $(n-1)$-dimensional volume $n \omega_n$. The integral average of an integrable function $f \colon X \to V$ over a measure $\mu$ and a $\mu$-measurable set $E \subset X$ with $\mu(E)>0$ is denoted by $\aint{ E } f \,d \mu$. We mainly apply this when $\mu = dz$ is the Lebesgue measure or the $(n-1)$-dimensional volume $\mu = ds$ of a sphere $\S^{n-1}(x,r)$.

Since we are working with maps $\R^n \to \R^m$, we identify linear maps $\R^n \to \R^m$ with matrices $\R^{ m \times n }$ through the standard basis. For a field of matrices $M \colon E \to \R^{ m \times n }$, the matrix at $v \in E$ is denoted by $M(v)$ and the evaluation of the matrix at a vector $w \in \R^n$ is denoted by $M(v)[w]$. When $\S \subset \R^n$, the definition is naturally extended to a map $M \colon \S \to \R^{ m \times (n-1) }$.

We denote by $g_{euc}$ the standard Euclidean inner product $g_{euc}(v,w) = \langle v, w \rangle$.

We also always assume that $3 \leq n \leq m$ and that $\omega$ denotes a calibration in $\Lambda^{n} \R^m$. Moreover, the notation $F \colon \R^n \supset \Omega \to ( \R^m, \omega )$ is reserved for conformal curves with only a few exceptions.

We also denote by $W^{k,p}_{loc}( \Omega, \R^m )$ and $W^{k,p}_{loc}( \Omega, \R^{ m \times n } )$ the standard distributional Sobolev spaces; the definition is reduced to the real-valued case by considering the standard bases. We only work with integers $k \in \N$ and $p \geq 1$.

\subsection{The natural differential}\label{section-natural-differential}
The theory of $n$-harmonic mappings is based on the natural differential $V_F = |DF|^{ \frac{n-2}{2} } DF$. In fact, $V_F \in W^{1,2}_{loc}( \Omega, \R^{ m \times n } )$, see e.g. \cite{uhlenbeck-1977-regularity-for-a-class-of-non-linear-elliptic-systems}. More precisely, we have the following.
\begin{proposition}[Caccioppoli's inequality]\label{proposition-natural-differential-caccioppoli}
    Let $3 \leq n \leq m$. Then there exists a constant $C = C(n) >0$ such that for every $n$-harmonic mapping $F \colon \R^n \supset \Omega \to \R^m$, it holds that $V_F \in W^{1,2}_{loc}( \Omega, \R^{ m \times n } )$ and
    \begin{align*}
        \int_{ \Omega } \eta^{2} |DV_F|^2 \,dz
        \leq
        C
        \inf_{ c \in \R^{ m \times n } }
        \int_{ \Omega } |D\eta|^2 \left| V_F - c \right|^2 \,dz
    \end{align*}
    for every $\eta \in C^{\infty}_{c}( \Omega )$.
\end{proposition}
\begin{proof}
    The argument is a standard application of the difference quotient method. First, the map
    \begin{align*}
        V \colon \R^{ m \times n } \to \R^{ m \times n },
        \quad M \mapsto |M|^{ \frac{n-2}{2} } M
    \end{align*}
    is a homeomorphism. In particular, if $c \in \R^{ m \times n }$, there exists a unique $M \in \R^{ m \times n }$ with $V(M) = c$. We also consider the homeomorphism
    \begin{align*}
        T \colon \R^{ m \times n } \to \R^{ m \times n },
        \quad M \mapsto |M|^{ n-2 } M.
    \end{align*}
    To complete the proof, we fix $\eta \in C^{\infty}_{c}(\Omega)$ and $v \in \S^{n-1}$ and test the $n$-harmonic equation of $F$ building upon the map
    \begin{align*}
        x \mapsto G(x) = F(x) - M[x]
    \end{align*}
    and $\eta$. To this end, for small enough $\mu > 0$, the support of $\eta$ is contained in
    \begin{align*}
        \Omega_{\mu} = \{ x \in \Omega \colon r(x) > \mu \}.
    \end{align*}
    For each number $0 < |h| < \mu$, we consider the map
    \begin{align*}
        \tau_{h,v} \colon \Omega_{\mu} \to \Omega, \quad x \mapsto x+ h v.
    \end{align*}
    We use the notation $g_{v,h}(x) = \frac{ g(x+hv) - g(x) }{ h }$ in the sequel.
    
    With this setup, we define the compactly-supported test function
    \begin{align*}
        \varphi = \eta^2 G_{v,h} \in W^{1,n}_{0}( \Omega, \R^m ).
    \end{align*}
    We use the commutativity of the difference quotient with the usual differentiation and the difference-quotient integration-by-parts to obtain the following equation
    \begin{align} \notag
        0 
        = \int_{\Omega} \langle D( \varphi_{-h,v} ), T(DF) \rangle \,dz 
        &= \int_{\Omega} \langle (D\varphi)_{-h,v}, T(DF) \rangle \,dz
        \\ \label{equation-integration-by-parts}
        &= \int_{\Omega} \langle D\varphi, ( T(DF) )_{h,v} \rangle \,dz
    \end{align}
    from the $n$-harmonicity. Here
    \begin{align*}
        D\varphi
        =
        2 \eta D\eta \otimes ( F_{h,v} - M[v] )
        +
        \eta^{2} ( DF )_{h,v}.
    \end{align*}
    Substituting to \eqref{equation-integration-by-parts} and multiplying both sides by $h^2$ gives
    \begin{gather*}
        \int_{ \Omega } \eta^{2} \langle DF \circ \tau_{h,v} - DF, ( T(DF) ) \circ \tau_{h,v} - T( DF ) \rangle \,dz
        \\
        =
        -2\int_{\Omega} \eta \langle D\eta \otimes ( F \circ \tau_{h,v} - F - h M[v] ), T( DF ) \circ \tau_{h,v} - T( DF ) \rangle \,dz.
    \end{gather*}
    We recall the standard monotonicity formulas for the $n$-harmonic system (see e.g. \cite[Section~10]{lindqvist-2006-notes-on-the-p-laplace-equation}): for a dimensional constant $c = c(n)$, it holds that
    \begin{align*}
        ( |A| + |B| )^{n-2} |A-B|^2
        &\leq
        c
        \langle T(A) - T(B), A-B \rangle
    \end{align*}
    and
    \begin{align*}
        \frac{ 1 }{ c } 
        ( |A| + |B| )^{ n-2 } 
        |A-B|^2
        \leq
        | V(A) - V(B) |^2
        \leq
        c
        ( |A| + |B| )^{ n-2 }
        |A-B|^2.
    \end{align*}
    Similarly,
    \begin{align*}
        | T(A) - T(B) | \leq ( |A| + |B| )^{n-2} |A-B|.
    \end{align*}
    We have, therefore, that
    \begin{gather*}
        \eta^{2} | V_F \circ \tau_{h,v} -  V_F |^2
        \leq
        c^2 
        \eta^{2}
        \langle T( DF ) \circ \tau_{h,v} - T( DF ), (DF) \circ \tau_{h,v} - DF \rangle.
    \end{gather*}
    Furthermore, since $\langle A \otimes B, C \rangle \leq |A||B||C|$, we have
    \begin{gather*}
        -2 \eta \langle D\eta \otimes ( F \circ \tau_{h,v} - F - h M[v] ), T( DF ) \circ \tau_{h,v} - T( DF ) \rangle
        \\
        \leq
        2 |\eta| |D\eta| | F \circ \tau_{h,v} - F - hM[v] | | T(DF) \circ \tau_{h,v} - T(DF) |
        \\
        \leq
        2 |\eta| |D\eta| | F \circ \tau_{h,v} - F - hM[v] | ( |DF| \circ \tau_{h,v} + |DF| )^{n-2} | DF \circ \tau_{h,v} - DF |. 
    \end{gather*}
    Hence,
    \begin{gather*}
        \int_{ \Omega }
            \eta^{2} | V_F \circ \tau_{h,v} -  V_F |^2 \,dz
        \\
        \leq
        2 c^2
        \int_{ \Omega }
            |\eta|| D\eta | | F \circ \tau_{h,v} - F - hM[v] | | T( DF )_{h,v} | \,dz.
    \end{gather*}
    By Young's inequality,
    \begin{align*}
        |\eta| |D\eta| | F \circ \tau_{h,v} - F - hM[v] | | DF \circ \tau_{h,v} - DF |
        &\leq
        \varepsilon \eta^2  | DF \circ \tau_{h,v} - DF |^2 
        \\
        &+
        C_\varepsilon |D\eta|^2 | F \circ \tau_{h,v} - F - hM[v] |^2.
    \end{align*}
    We note that
    \begin{gather*}
        ( |DF| \circ \tau_{h,v} + |DF| )^{n-2}| DF \circ \tau_{h,v} - DF |^2
        \leq
        c |V_F \circ \tau_{h,v} - V_F|^2.
    \end{gather*}
    Therefore, for small enough $\varepsilon = \varepsilon(n)$, we conclude
    \begin{gather*}
        \int_{ \Omega }
            \eta^{2} | V_F \circ \tau_{h,v} -  V_F |^2 \,dz
        \\
        \leq
        \frac{ 2 c^2 C_{\varepsilon} }{ 1 - 2 \varepsilon c^3 }
        \int_{\Omega} |D\eta|^2 ( |DF| \circ \tau_{h,v} + |DF| )^{n-2}| F \circ \tau_{h,v} - F - hM[v] |^2 \,dz.
    \end{gather*}
    We divide by $h^2$, pass to the limit $h \to 0$, and use the difference quotient characterization of Sobolev spaces to deduce
    \begin{align*}
        \int_{ \Omega }
            \eta^{2} | DV_F[v] |^2 \,dz
        \leq
        \frac{ 2^{n-1} c^2 C_{\varepsilon} }{ 1 - 2 \varepsilon c^3 }
        \int_{\Omega} |D\eta|^2 |DF|^{n-2}| DF[v] - M[v] |^2 \,dz.
    \end{align*}
    We consider $v = e_1,\dots,e_n$ and sum the results to obtain
    \begin{align*}
        \int_{ \Omega }
            \eta^{2} | DV_F |^2 \,dz
        &\leq
        \frac{ 2^{n-1} c^2 C_{\varepsilon} }{ 1 - 2 \varepsilon c^3 }
        \int_{\Omega} |D\eta|^2 |DF|^{n-2}| DF - M |^2 \,dz
        \\
        &\leq
        \frac{ 2^{n-1} c^2 C_{\varepsilon} }{ 1 - 2 \varepsilon c^3 }
        \int_{\Omega} |D\eta|^2 ( |DF| + |M| )^{n-2}| DF - M |^2 \,dz
        \\
        &\leq
        \frac{ 2^{n-1} c^3 C_{\varepsilon} }{ 1 - 2 \varepsilon c^3 }
        \int_{\Omega} |D\eta|^2 | V_F - c |^2 \,dz.
    \end{align*}
    The claim follows.
\end{proof}
We have the following standard corollary.
\begin{corollary}\label{corollary-improved-integrability}
    Let $3 \leq n \leq m$. Then there exists $\delta = \delta(n) > 0$ such that $|DV_F| \in L^{ 2(1+\delta) }_{loc}( \Omega )$ for every $n$-harmonic mapping $F \colon \R^n \supset \Omega \to \R^m$.
\end{corollary}
\begin{proof}
    We apply Proposition~\ref{proposition-natural-differential-caccioppoli} with a standard cut-off function $\phi \in C^{\infty}_{c}( \Omega )$ with range $[0,1]$, zero outside $B(x,2r)$, one in $B(x,r)$, and differential bound $|D\phi| \leq 2/r$. Then the proposition implies that
    \begin{align*}
        \int_{ B(x,r) } |DV_F|^2 \,dz
        \leq
        \frac{ 4 C }{ r^2 }
        \int_{ B(x,2r) } \left| V_F - \aint{ B(x,2r) } V_F \,dz \right|^2 \,dz.
    \end{align*}
    Dividing both sides by the measure of $B(x,r)$, we obtain
    \begin{align*}
        \left( \aint{ B(x,r) } |DV_F|^2 \,dz \right)^{ \frac{1}{2} }
        \leq
        \left( \frac{2}{r} \right)
        \sqrt{ 2^{n}C }
        \left( \aint{ B(x,2r) } \left| V_F - \aint{ B(x,2r) } V_F \,dz \right|^2 \,dz \right)^{ \frac{1}{2} }.
    \end{align*}
    Applying the local Sobolev--Poincaré inequality for the exponent $p = 2n/(n+2) < 2$ to $V_F \in W^{1,2}_{loc}( \Omega, \R^{m \times n} ) \subset W^{1,p}_{loc}( \Omega, \R^{ m \times n } )$ implies that
    \begin{align*}
        \left( \aint{ B(x,2r) } \left| V_F - \aint{ B(x,2r) } V_F \,dz \right|^2 \,dz \right)^{ \frac{1}{2} }
        \leq
        d r
        \left( \aint{ B(x,2r) } | DV_F |^{p} \,dz \right)^{ \frac{1}{p} }
    \end{align*}
    for a constant $d$ depending only on $n$. By combining the inequalities, it follows that
    \begin{align*}
        \left( \aint{ B(x,r) } |DV_F|^2 \,dz \right)^{ \frac{1}{2} }
        \leq
        \left( 2d \right)
        \sqrt{ 2^{n}C }
        \left( \aint{ B(x,2r) } | DV_F |^{p} \,dz \right)^{ \frac{1}{p} }.
    \end{align*}
    Now, the standard Gehring lemma implies the claim, see~\cite{gehring-1973-the-lp-integrability-of-the-partial-derivatives-of-a-quasiconformal-mapping,iwaniec-1998-the-gehring-lemma}.
\end{proof}

\section{The higher-order regularity of the natural differential}\label{section-higher-order-regularity}
In this section, we prove the following theorems. First, a second-order result for the function $|V_F|^2$.
\begin{theorem}\label{theorem-higher-order-regularity}
    Let $3 \leq n \leq m$ and let $\delta = \delta(n)$ be as in Corollary~\ref{corollary-improved-integrability}. Then $|V_F|^2 \in W^{2,1+\delta}_{loc}( \Omega )$ for every conformal curve $F \colon \R^n \supset \Omega \to ( \R^m, \omega )$.
\end{theorem}
The second theorem is closely related to Theorem~\ref{theorem-higher-order-regularity}. Formulating the result requires further setup. We first consider the matrix $Q^F = (Q_{ji}^F)$ defined by the identity
\begin{align*}
    Q_{ji}^F 
    = 
    \langle \partial_{j}V_F, \partial_{i}V_F \rangle 
    - 
    \theta^2 \partial_{j}|V_F|\partial_{i}|V_F|
\end{align*}
for the dimensional constant
\begin{equation}\label{equation-dimensional-constant}
    \theta = \frac{n-2}{n}.
\end{equation}
We later need the related constant
\begin{align}\label{equation-dimensional-constant-two}
    k = 1-\theta^2 = 4 \frac{n-1}{n^2}.
\end{align}
With this setup, we define the trace-adjusted matrix field
\begin{align*}
    \mathcal{T}^{F} = Q^F - \frac{1}{2} \mathrm{trace}(Q^F) I_{n \times n},
\end{align*}
where $I_{ n \times n }$ is the identity matrix.
\begin{theorem}\label{theorem-stress-energy-tensor}
    Let $3 \leq n \leq m$ and consider a conformal curve $F \colon \R^n \supset \Omega \to ( \R^m, \omega )$. Then
    \begin{align*}
        Q^F 
        = 
        \frac{n-2}{n^2} D^{2} |V_F|^2 + \frac{ \Delta( |V_F|^2 ) }{ n^2 } I_{ n \times n }
    \end{align*}
    and
    \begin{align}\label{equation-stress-energy-tensor-two-representations}
        \mathcal{T}^F 
        = 
        \frac{n-2}{n^2}\left( D^{2}( |V_F|^2 ) - \Delta( |V_F|^2 ) I_{ n \times n } \right)
    \end{align}
    everywhere in the full-rank set $U_F$ and almost everywhere on the singular set $S_F$. Furthermore, $\mathcal{T}^F$ is distributionally divergence-free. That is, $\sum_{j} \partial_{j}( \mathcal{T}_{ij}^F ) = 0$ distributionally for every $i$.
\end{theorem}
\begin{corollary}\label{corollary-representation-of-the-laplacian}
    It holds that
    \begin{align*}
        |DV_F|^2
        -
        \theta^2
        |D|V_F||^2
        =
        \frac{k}{2} \Delta( |V_F|^2 )
    \end{align*}
    almost everywhere on $\Omega$ for the dimensional constants from \eqref{equation-dimensional-constant} and \eqref{equation-dimensional-constant-two}.
\end{corollary}
\begin{proof}
    The claim is immediate by taking a trace of the two representations of the matrix field $Q^F$ (or from \eqref{equation-stress-energy-tensor-two-representations}).
\end{proof}

The proof of the theorems builds upon the following representation result for the distributional Laplacian of $|V_F|^2$.
\begin{proposition}\label{proposition-initial-information-on-the-subharmonicity}
    If $F \colon \R^n \supset \Omega \to ( \R^m, \omega )$ is a conformal curve for $3 \leq n \leq m$, then $|DF|^{ \frac{n-2}{2} } $ is subharmonic and in $W^{1,2}_{loc}( \Omega )$. In particular, $|V_F|^2$ is subharmonic and its distributional Laplacian is absolutely continuous with respect to the Lebesgue measure with a density vanishing almost everywhere on $S_F$.
\end{proposition}
\begin{proof}
    The subharmonicity of $\rho_F = |DF|^{ \frac{n-2}{2} }$ is proved in \cite[Theorem 1.1]{ikonen-pankka-2024-liouvilles-theorem-in-calibrated-geometries} based on the isoperimetric inequality for conformal curves \cite[Theorem~4.1]{ikonen-pankka-2024-liouvilles-theorem-in-calibrated-geometries} (see also \cite{almgren-1986-optimal-isoperimetric-inequalities,ikonen-2026-pushforward-of-currents-under-sobolev-maps}). Since $\rho_F$ is continuous, Caccioppoli's inequality for subharmonic functions gives that $\rho_F \in W^{1,2}_{loc}( \Omega )$. Moreover, since $|V_F|^2 = \rho_{F}^{ \frac{2n}{n-2} }$ (where $t \mapsto t^{ \frac{2n}{n-2} }$ is $C^2$-regular and convex in $[0,\infty)$), the representation
    \begin{align*}    
        \Delta( |V_F|^2 )
        = 
        \frac{2n}{n-2}
        \left( 
            \rho^{ \frac{n+2}{n-2} }_F \Delta\rho_F
            +
            \frac{n+2}{n-2}
            \rho^{ \frac{4}{n-2} }_F |D\rho_F|^2 \,dz
        \right)
    \end{align*}
    follows from the chain rule; the equality should be understood as Radon measures. Now, since $F$ is real-analytic in the full-rank set $U_F$, it follows that $\Delta( |V_F|^2 )$ is absolutely continuous with respect to the Lebesgue measure on $U_F$. Also, the measure $\rho^{ \frac{n+2}{n-2} }_F \Delta\rho_F$ vanishes on $S_F$ since $S_F = \{ \rho_F = 0 \}$. Here $\rho^{ \frac{4}{n-2} }_F |D\rho_F|^2 \,dz$ gives zero measure to $S_F$. Hence, $\Delta( |V_F|^2 )$ vanishes on $S_F$. The absolute continuity follows.
\end{proof}

Proposition~\ref{proposition-initial-information-on-the-subharmonicity} shows that the Laplacian $\Delta( |V_F|^2 )$ is completely determined by the full-rank set $U_F$. The remainder of the section is devoted to establishing Theorem~\ref{theorem-stress-energy-tensor} in $U_F$, formulated as follows.
\begin{proposition}\label{proposition-variance-of-V}
    In the full-rank set $U_F$, a conformal curve satisfies
    \begin{align*}
        Q^F 
        = 
        \frac{n-2}{n^2} D^{2} |V_F|^2 + \frac{ \Delta( |V_F|^2 ) }{ n^2 } I_{ n \times n }.
    \end{align*}
    In particular,
    \begin{align*}
        \mathcal{T}^F 
        = 
        \frac{n-2}{n^2}\left( D^{2}( |V_F|^2 ) - \Delta( |V_F|^2 ) I_{ n \times n } \right).
    \end{align*}
\end{proposition}
    To highlight the importance of Propositions~\ref{proposition-initial-information-on-the-subharmonicity} and \ref{proposition-variance-of-V}, we present a proof of Theorems~\ref{theorem-higher-order-regularity} and \ref{theorem-stress-energy-tensor} using the proposition.
\begin{proof}[Proof of Theorem~\ref{theorem-higher-order-regularity} assuming Proposition~\ref{proposition-variance-of-V}]
    Taking the trace of $Q^F$ using the representation of Proposition~\ref{proposition-variance-of-V} leads to the identity
    \begin{align}\label{equation-laplacian-and-energy}
        |DV_F|^2 - \theta^2 |D|V_F||^2
        =
        2 \frac{n-1}{n^2} \Delta( |V_F|^2 )
    \end{align}
    on $U_F$. Since $V_F \in W^{1,2}_{loc}( \Omega, \R^{ m \times n } )$ by Proposition~\ref{proposition-natural-differential-caccioppoli}, the strong locality of the Sobolev differential leads to $|DV_F| = 0 = |D|V_F||$ almost everywhere on $S_F$. Since by Proposition~\ref{proposition-initial-information-on-the-subharmonicity}, it holds that $\Delta( |V_F|^2 ) = 0$ almost everywhere on $S_F$, we conclude that \eqref{equation-laplacian-and-energy} holds almost everywhere on $\Omega$. Since the left-hand side of \eqref{equation-laplacian-and-energy} is comparable to $|DV_F|^2$ and $|DV_F|^2 \in L^{1+\delta}_{loc}( \Omega )$ by Corollary~\ref{corollary-improved-integrability}, we see that $\Delta( |V_F|^2 ) \in L^{1+\delta}_{loc}(\Omega)$. Recall that, for any ball $\overline{B}(x,r) \subset \Omega$,
    \begin{align*}
        |V_F|^2 = h + s \quad\text{in $B(x,r)$}
    \end{align*}
    for harmonic $h \colon B(x,r) \to \R$ and subharmonic $s$, where $s$ is the Newtonian potential of $\Delta( |V_F|^2 ) \in L^{1+\delta}( B(x,r) )$. Since $h$ is smooth, the conclusion $|V_F|^2 \in W^{2,1+\delta}( B(x,r) )$ follows from the regularity theory of the Poisson equation, see e.g. \cite[Chapter~9, Theorem~9.9]{gilbarg-trudinger-2001-elliptic-partial-differential-equations-of-second-order}. The claim $|V_F|^2 \in W^{2,1+\delta}_{loc}( \Omega )$ follows by a covering argument.
\end{proof}
\begin{proof}[Proof of Theorem~\ref{theorem-stress-energy-tensor} assuming Proposition~\ref{proposition-variance-of-V}]
    By Theorem~\ref{theorem-higher-order-regularity} and strong locality of the Sobolev differential, it holds
    \begin{align*}
        \frac{n-2}{n^2} D^{2} |V_F|^2 + \frac{ \Delta( |V_F|^2 ) }{ n^2 } I_{ n \times n } = 0
    \end{align*}
    almost everywhere on $S_F$. Similarly, $Q^{F} = 0$ almost everywhere on $S_F$ by $V_F \in W^{1,2}_{loc}( \Omega, \R^{m\times n} )$. Hence the claimed identity for $Q^F$ follows from Proposition~\ref{proposition-variance-of-V} which immediately also yields the corresponding identity for $\mathcal{T}^F$.

    It remains to verify that $\mathcal{T}^F$ is distributionally divergence-free. To this end, notice that the distributional derivatives of $\mathcal{T}^F$ are insensitive to changing $\mathcal{T}^F$ in a set negligible for the Lebesgue measure. Hence, using the representation \eqref{equation-stress-energy-tensor-two-representations} and the commutativity and linearity of distributional partial derivatives, it holds that
    \begin{align*}
        \sum_{j} \partial_{j}( \mathcal{T}^F_{ij} )
        &=
        \frac{n-2}{n^2}
        \sum_{j}
        \left(
            \partial_{i}\partial_{j}\partial_{j}(|V_F|^2) - \partial_{j}( \Delta(|V_F|^2) )\delta_{ij}
        \right)
        \\
        &=
        \frac{n-2}{n^2}
        \left(
            \partial_{i} \sum_{j}\partial_{j}\partial_{j}(|V_F|^2) - \partial_{i}( \Delta(|V_F|^2) )
        \right)
        =
        0.
    \end{align*}
    The conclusion follows.
\end{proof}

    The remaining section is dedicated to proving Proposition~\ref{proposition-variance-of-V}. This is achieved in a series of lemmas. To set up the statements, we let $s = |DF|$, $f = \log s$, $b = (n-2)/2$, $V = s^{b} DF$ (noting $V = V_F$), and $w = s^n$.
\begin{lemma}\label{lemma-conformal-curve-re-expression}
    In the full-rank set $U_F$, it holds that
    \begin{align}\label{equation-conformality-and-n-harmonicity}
        \langle \partial_j F, \partial_{i} F \rangle = \frac{ e^{2f} }{ n } \delta_{ji}
        \quad\text{and}\quad
        \Delta F
        =
        -
        (n-2) \sum_{l} (\partial_{l}f) (\partial_{l}F ),
    \end{align}
    where $\delta_{ji}$ is the Kronecker delta.
\end{lemma}
\begin{proof}
    The conformality of $F$ is immediate from the definition (see e.g. \cite[Lemma~2.1]{ikonen-pankka-2024-liouvilles-theorem-in-calibrated-geometries}) and the distributional $n$-harmonicity is proved for instance in \cite[Theorem~5.1]{heikkila-pankka-prywes-2023-quasiregular-curves-of-small-distortion-in-product-manifolds}. In the full-rank set, the $n$-harmonic equation $\mathrm{div}( s^{n-2} DF ) = 0$ can be expressed as stated using the product rule.
\end{proof}
    Differentiation leads to the following lemma.
\begin{lemma}\label{lemma-conformal-curve-3rd-order-n-harmonicty}
    In the full-rank set $U_F$, it holds that
    \begin{align} \label{equation-laplacian-formula}
        \langle \partial_{j} \Delta F, \partial_{i}F \rangle 
        &= 
        - 
        (n-2)
        \frac{ e^{2f} }{ n }
        \left(
            (\partial_{j}\partial_{i}f)
            +
            |Df|^2 \delta_{ji}
        \right)
    \end{align}
    and that
    \begin{align}\label{equation-cyclic-differentiation}
        \langle \partial_{j}\partial_{l}F, \partial_{i}F \rangle
        =
        \frac{1}{n}e^{2f}
        \left(
            (\partial_{l}f) \delta_{ji}
            +
            (\partial_{j}f) \delta_{il}
            -
            (\partial_{i}f) \delta_{lj}
        \right).
    \end{align}
\end{lemma}
\begin{proof}
    Regarding the second identity, we differentiate the conformality equation, in the formulation of \eqref{equation-conformality-and-n-harmonicity}, to get
    \begin{align*}
        \langle \partial_{j}\partial_{l}F, \partial_{i}F \rangle
        +
        \langle \partial_{j}F, \partial_{i} \partial_{l}F \rangle
        =
        \frac{2}{n}
        e^{2f} (\partial_{l}f) \delta_{ji}.
    \end{align*}
    Cyclically permuting the indices and summing leads to the assertion.
    
    For the first identity, we differentiate the $n$-harmonic equation and apply the conformality to obtain
    \begin{align*}
        \langle \partial_{j} \Delta F, \partial_{i}F \rangle 
        &=
        - (n-2)
        \left\langle 
            \sum_{l} \partial_{j}( (\partial_{l}f) (\partial_{l}F ) ), 
            \partial_{i}F 
        \right\rangle
        \\ \notag
        &=
        - 
        (n-2)
        \frac{ e^{2f} }{ n }
        (\partial_{j}\partial_{i}f)
        -
        (n-2)
        \sum_{l}
        \langle (\partial_{l}f)( \partial_{j}\partial_{l}F ), \partial_{i}F \rangle
        \\
        &=
        - 
        (n-2)
        \frac{ e^{2f} }{ n }
        (\partial_{j}\partial_{i}f)
        -
        (n-2)
        \frac{ e^{2f} }{ n }
        |Df|^{2}\delta_{ji}.
    \end{align*}
\end{proof}
    We further differentiate \eqref{equation-cyclic-differentiation} and obtain the following higher-order equation.
\begin{lemma}\label{lemma-second-order-differential-identity}
    In the full-rank set $U_F$, it holds that
    \begin{align*}
        \sum_{l}
        \langle \partial_{j}\partial_{l}F, \partial_{i}\partial_{l}F \rangle
        &=
        \frac{ e^{2f} }{ n }
        \left(
            (n-2) (\partial_{j}\partial_{i}f)
            +
            \left(
                n |Df|^2
                +
                \Delta f
            \right)
            \delta_{ji}
        \right).
    \end{align*}
\end{lemma}
\begin{proof}
Differentiating \eqref{equation-cyclic-differentiation} with $\partial_{l}$ and summing leads to
    \begin{align*}
        \langle \partial_{j}\Delta F, \partial_{i}F \rangle
        +
        \sum_{l}
        \langle \partial_{j}\partial_{l}F, \partial_{i}\partial_{l}F \rangle
        &=
        \sum_{l} 
        \frac{ e^{2f} }{ n } 2(\partial_{l}f)
        \left(
            (\partial_{l}f) \delta_{ji}
            +
            (\partial_{j}f) \delta_{il}
            -
            (\partial_{i}f) \delta_{lj}
        \right)
        \\
        &+
        \sum_{l}
        \frac{ e^{2f} }{ n }
        \left(
            (\partial_{l}\partial_{l}f) \delta_{ji}
            +
            (\partial_{j}\partial_{l}f) \delta_{il}
            -
            (\partial_{i}\partial_{l}f) \delta_{lj}
        \right)
        \\
        &=
        \frac{ e^{2f} }{ n } 2 |Df|^2 \delta_{ji}
        +
        \frac{ e^{2f} }{ n }
        \Delta f \delta_{ji}.
    \end{align*}
    Using \eqref{equation-laplacian-formula} leads to the claim.
\end{proof}

\begin{proof}[Proof of Proposition~\ref{proposition-variance-of-V}]
    Recall $V = s^{b} DF = e^{bf} DF$. The product rule gives
    \begin{align*}
        \partial_{j}V = e^{bf}( b(\partial_{j}f)DF + D( \partial_{j}F ) )
    \end{align*}
    so
    \begin{align*}
        \partial_{j}V_{k}[e_l] = e^{bf}( b(\partial_{j}f) (\partial_l F_k) + (\partial_{l}\partial_{j}F_k) ).
    \end{align*}
    Hence, recalling $2b = n-2$,
    \begin{align*}
        \langle \partial_{j}V, \partial_{i}V \rangle
        &=
        e^{(n-2)f}
        \sum_{k,l} 
            \left(
                b(\partial_{j}f) (\partial_{l} F_k) 
                + 
                (\partial_{l}\partial_{j}F_k)
            \right)
            \left(
                b(\partial_{i}f) (\partial_l F_k) 
                + 
                (\partial_{l}\partial_{i}F_k)
            \right)
        \\
        &=
        e^{(n-2)f}
        \left( 
            b^{2}(\partial_{j}f)(\partial_{i}f)|DF|^2
            +
            \sum_{l}
            \langle \partial_{l}\partial_{j}F, \partial_{l}\partial_{i}F \rangle
        \right)
        \\
        &
        +
        b
        e^{(n-2)f}
        \left(
            (\partial_{j}f)
            \sum_{k,l} (\partial_{l} F_k)(\partial_{l}\partial_{i}F_k)
            +
            (\partial_{i}f)
            \sum_{k,l} (\partial_{l} F_k)(\partial_{l}\partial_{j}F_k)
        \right).
    \end{align*}
    Recall that $|DF|^2 = e^{2f}$. We simplify $\sum_{l}\langle \partial_{l}\partial_{j}F, \partial_{l}\partial_{i}F \rangle$ using Lemma~\ref{lemma-second-order-differential-identity} and use \eqref{equation-cyclic-differentiation} to simplify $\sum_{k,l} (\partial_{l} F_k)(\partial_{l}\partial_{i}F_k)$ and $\sum_{k,l} (\partial_{l} F_k)(\partial_{l}\partial_{j}F_k)$. These lead to
    \begin{align*}
        \langle \partial_{j}V, \partial_{i}V \rangle
        &=
        e^{(n-2)f}
        \left( 
            b^2(\partial_{j}f)(\partial_{i}f)e^{2f}
            +
            \frac{ e^{2f} }{ n }
            \left(
                (n-2) (\partial_{j}\partial_{i}f)
                +
                \left(
                    n |Df|^2
                    +
                    \Delta f
                \right)
                \delta_{ji}
        \right)
        \right)
        \\
        &+
        b
        e^{(n-2)f}
        \left(
            (\partial_{j}f)
            e^{2f}
            \left(
                \partial_{i}f
            \right)
            +
            (\partial_{i}f)
            e^{2f}
            \left(
                \partial_{j}f
            \right)
        \right)
        \\
        &=
        e^{nf}
        \left(
            (b^2+2b)(\partial_{j}f)(\partial_{i}f)
            +
            \frac{n-2}{n} (\partial_{j}\partial_{i}f)
            +
            \left(
                |Df|^2
                +
                \frac{ \Delta f }{ n }
            \right)
            \delta_{ji}
        \right)
        \\
        &=
        e^{nf}
        \left(
            \frac{n^2-4}{4}(\partial_{j}f)(\partial_{i}f)
            +
            \frac{n-2}{n} (\partial_{j}\partial_{i}f)
            +
            \left(
                |Df|^2
                +
                \frac{ \Delta f }{ n }
            \right)
            \delta_{ji}
        \right).
    \end{align*}
    Here $f = \log(s) = \frac{1}{n} \log(w)$ so
    \begin{align*}
        \partial_{j}f = \frac{1}{n}\frac{ \partial_{j}w }{ w }
        \quad\text{and}\quad
        \partial_{j}\partial_{i}f = -\frac{1}{n} \frac{ (\partial_{j}w)(\partial_{i}w) }{ w^2 } + \frac{1}{n} \frac{ \partial_{j}\partial_{i}w }{ w }.
    \end{align*}
    Hence,
    \begin{align*}
        \langle \partial_{j}V, \partial_{i}V \rangle
        &=
        \frac{w}{n^2}
        \left(
            \frac{n^2-4}{4}
            \frac{ (\partial_{j}w)(\partial_{i}w) }{ w^2 }
            -
            (n-2) \frac{ (\partial_{j}w)(\partial_{i}w) }{ w^2 }
            +
            (n-2) \frac{ \partial_{j}\partial_{i}w }{ w }
        \right)
        \\
        &+
        \frac{w}{n^2}
        \left(
            \frac{ |Dw|^2 }{ w^2 }
            -
            \frac{ |Dw|^2 }{ w^2 }
            +
            \frac{ \Delta w }{ w }
        \right)
        \delta_{ji}
        \\
        &=
        b^2
        \frac{ (\partial_{j}w)(\partial_{i}w) }{ n^2 w }
        +
        \frac{n-2}{n^2} \partial_{j}\partial_{i}w
        +
        \frac{ \Delta w }{ n^2 }
        \delta_{ji}.
    \end{align*}
    Recall the definitions $\theta = (n-2)/n$ and $w = |V|^2$. Then $(\partial_{j}w)(\partial_{i}w) = (4w) (\partial_{j}|V|)(\partial_{i}|V|)$ leads to
    \begin{align*}
        \langle \partial_{j}V, \partial_{i}V \rangle
        &=
        \theta^{2} (\partial_{j}|V|)(\partial_{i}|V|)
        +
        \frac{n-2}{n^2} \partial_{j}\partial_{i}w
        +
        \frac{ \Delta w }{ n^2 }
        \delta_{ji}.
    \end{align*}
    That is,
    \begin{align*}
        Q_{ji}
        &=
        \frac{n-2}{n^2} \partial_{j}\partial_{i}w
        +
        \frac{ \Delta w }{ n^2 }
        \delta_{ji}
    \end{align*}
    and $$\mathrm{trace}(Q) = \left( \frac{n-2}{n^2} \right) \Delta w + \left( \frac{n}{n^2} \right) \Delta w = \frac{ 2(n-1) }{ n^2 } \Delta w.$$ Hence,
    \begin{align*}
        \mathcal{T}^F_{ji} 
        &= 
        Q_{ji} - \frac{1}{2} \mathrm{trace}(Q) \delta_{ji}
        \\
        &=         
        \frac{n-2}{n^2} \partial_{j}\partial_{i}w
        +
        \frac{ \Delta w }{ n^2 }
        \delta_{ji}
        -
        \frac{n-1}{n^2} \Delta w
        \delta_{ji}
        \\
        &=
        \frac{n-2}{n^2}\left( \partial_{j}\partial_{i} w - \Delta w \delta_{ji} \right).
    \end{align*}
    The claim follows.
    \end{proof}

\section{Frequency function for conformal curves}\label{section-frequency-function-for-conformal-curves}
    For this section, we fix a conformal curve $F \colon \R^n \supset \Omega \to ( \R^m, \omega )$ for a calibration $\omega \in \Lambda^{n} \R^m$ and $3 \leq n \leq m$. We derive important consequences of Theorem~\ref{theorem-stress-energy-tensor}. More concretely, we prove the monotonicity formula, Theorem~\ref{theorem-frequency-function-monotonicity} and the unique continuation principle, Theorem~\ref{theorem-unique-continuation-principle-open-set}.
    
    We recall the dimensional constants
    \begin{align*}
        \theta = \frac{n-2}{n}
        \quad\text{and}\quad
        k = 1 - \theta^2 = 4 \frac{n-1}{n^2}
    \end{align*}
    from \eqref{equation-dimensional-constant} and \eqref{equation-dimensional-constant-two}. We introduce additional functions that are used to analyze the singularities of $F$. The first is defined by the boundary energy
    \begin{equation*}
        H_{F}(x,r) = \int_{ \S^{n-1}(x,r) } |V_F|^2 \,ds = \int_{ \S^{n-1}(x,r) } |DF|^n \,ds
    \end{equation*}
    and its normalized version
    \begin{equation*}
        h_{F}(x,r)
        =
        \frac{ H_{F}(x,r) }{ n \omega_n r^{n-1} }.
    \end{equation*}
    The second function is the integral of the Laplacian of $|V_F|^2$ over balls. A natural normalization is as follows
    \begin{equation*}
        e_{F}(x,r)
        =
        \frac{r}{n} \aint{ B(x,r) } \Delta( |V_F|^2 ) \,dz.
    \end{equation*}
        Since $|V_F|^2$ is subharmonic, the Poisson--Jensen formula states that
        \begin{align}\label{equation-first-fundamental-theorem-of-nevanlinna-theory-prerequisite}
            h_{F}(x,r)
            =
            |V_{F}|^2(x)
            +
            \int_{ 0 }^{ r } e_{F}(x,a) \,da.
        \end{align}
    \begin{remark}\label{remark-monotonicity-formula-deconstrunction}
        By \eqref{equation-first-fundamental-theorem-of-nevanlinna-theory-prerequisite}, it follows that $r \mapsto h_F(x,r)$ is non-decreasing. Moreover, the Poisson--Jensen formula implies that each of the sets
        \begin{align*}
           & \{ h_{F}(x,\cdot) = 0 \} \coloneqq \{ r \in (0,r(x)) \colon h_{F}(x,r) = 0 \},
            \\
            &\{ h_{F}(x,\cdot) > 0 \} \coloneqq \{ r \in (0,r(x)) \colon h_{F}(x,r) > 0 \},
            \\
            &\{ h_{F}(x,\cdot) > 0 \} \cap \{ e_{F}(x,\cdot) = 0 \},
            \quad\text{and}
            \\
            &\{ h_{F}(x,\cdot) > 0 \} \cap \{ e_{F}(x,\cdot) > 0 \}
        \end{align*}
        is either empty or an interval, where $\{ e_{F}(x,\cdot) = 0 \}$ and $\{ e_{F}(x,\cdot) > 0 \}$ are defined in the obvious way. The representation \eqref{equation-first-fundamental-theorem-of-nevanlinna-theory-prerequisite} of $h_F$ yields
        \begin{align*}
            \{ h_{F}(x,\cdot) = 0 \} \subset \{ e_{F}(x,\cdot) = 0 \}.
        \end{align*}
    \end{remark}
    In the following sections, we need much finer analysis of \eqref{equation-first-fundamental-theorem-of-nevanlinna-theory-prerequisite} than the one in Remark~\ref{remark-monotonicity-formula-deconstrunction}. In fact, by definition of the frequency and \eqref{equation-first-fundamental-theorem-of-nevanlinna-theory-prerequisite}, we have
        \begin{align}\label{equation-first-fundamental-theorem-of-nevanlinna-theory}
            h_{F}(x,r)
            &=
            |V_F|^2(x)
            +
            n
            \int_{0}^{r} \frac{1}{a} \mathcal{N}_{F}(x,a) h_{F}(x,a) \,da.
        \end{align}
        It follows that $\mathcal{N}_{F}(x,r)$ is related to the convergence rate of $h_{F}(x,r)$ to $|V_F|^2(x)$ as $r \to 0^{+}$, justifying the term \emph{frequency}.
        
        The frequency interpretation is also clear from the following doubling property.
        \begin{theorem}\label{theorem-first-fundamental-theorem-of-nevanlinna-theory-and-frequency}
            Let $0<r_0<r(x)$. For every $0 < s < r < r_0$, it holds that
            \begin{align*}
                0
                \leq
                h_{F}(x,r)
                -
                h_{F}(x,s)
                &\leq
                n
                \mathcal{N}_{F}(x,r_0)
                \log\left( \frac{r}{s} \right)
                h_{F}(x,r).
            \end{align*}
        \end{theorem}
    \begin{proof}[Proof of Theorem~\ref{theorem-first-fundamental-theorem-of-nevanlinna-theory-and-frequency} assuming Theorem~\ref{theorem-frequency-function-monotonicity}] 
        We apply \eqref{equation-first-fundamental-theorem-of-nevanlinna-theory} and the monotonicity of both the frequency and $h_{F}$ to obtain the inequality.
    \end{proof}

\subsection{The growth of the auxiliary functions}\label{section-growth-of-the-auxiliary-functions}
    In this section, we analyze the functions $H_{F}(x,r)$ and
    \begin{align*}
        E_F(x,r)
        &=
        \int_{ B(x,r) }
            |DV_F|^2
            -
            \theta^2
            |D|V_F||^2
        \,dz.
    \end{align*}
    We note that Corollary~\ref{corollary-representation-of-the-laplacian} leads to the equivalent definition
    \begin{align}\label{equation-energy-in-terms-of-laplacian}
        E_F(x,r) = \frac{k}{2} \int_{ B(x,r) } \Delta( |V_F|^2 ) \,dz.
    \end{align}
    Using \eqref{equation-energy-in-terms-of-laplacian}, it holds that
    \begin{equation}\label{equation-alternate-expression-of-frequency}
        \frac{n}{2} \mathcal{N}_{F}(x,r) = \frac{ r E_{F}(x,r) }{ k H_{F}(x,r) }
    \end{equation}
    whenever $h_{F}(x,r)>0$ and zero otherwise. Hence, to analyze the monotonicity of the frequency, it suffices to differentiate the auxiliary functions $E_F$ and $H_F$.
    
    To simplify notation, we omit $x$ and $F$ from the computations below and we suppose $x = 0$. To recover the more general formulas, the radial derivatives $\partial_{r}$ need to be replaced by the directional derivative in the direction of the unit vector field $y \mapsto (y-x)/|y-x|$ for $y \neq x$.
    \begin{lemma}[A radial Rellich--Pohozaev type identity]\label{lemma-key-identities-for-the-energy}
        For almost every $r$, it holds that
        \begin{align*}
            E'(r)
            &=
            \frac{n-2}{r} E(r)
            +
            2
            \int_{ \S^{n-1}(x,r) }
                |\partial_{r} V|^2 - \theta^2 |\partial_{r}|V||^2
            \,ds.
        \end{align*}
    \end{lemma}
    \begin{proof}
        For the proof, we recall the assumption $x=0$. The claim follows by considering the tensor $\mathcal{T}^F$ from \eqref{equation-stress-energy-tensor-two-representations}. By construction,
        \begin{align*}
            \mathrm{trace}( \mathcal{T}^F ) 
            =
            -
            \frac{n-2}{2}
            \left(
                |DV|^2
                -
                \theta^2
                |D|V||^2
            \right).
        \end{align*}
        In particular, using Theorem~\ref{theorem-stress-energy-tensor}, the distributional divergence of the contraction $$y \mapsto \mathcal{T}^F(y)[y] = \sum_{i} \sum_{j} y_j \mathcal{T}^F_{ji}(y) e_i$$ is equal to $y \mapsto \mathrm{trace}( \mathcal{T}^F(y) )$. Indeed, as the distributional divergence of $\mathcal{T}^F$ is zero, the chain rule yields that the distributional divergence of $y \mapsto \mathcal{T}^F(y)[y]$ is equal to the claimed quantity. By a standard mollification argument for $\mathcal{T}^F$, the divergence theorem guarantees for almost every $r$ that
    \begin{align*}
        0
        &=
        \frac{n-2}{2}
        E(r)
        +
        \int_{ \S^{n-1}(x,r) }
            \left\langle \mathcal{T}^F(y)[y], \frac{y}{r} \right\rangle
        \,ds(y)
        \\
        &=
        \frac{n-2}{2} E(r)
        \\
        &+
        r
        \int_{ \S^{n-1}(x,r) }
            |\partial_{r} V|^2 - \theta^2 |\partial_{r}|V||^2
            -
            \frac{1}{2}
            \left(
                |DV|^2
                -
                \theta^2
                |D|V||^2
            \right)
        \,ds
        \\
        &=
        \frac{n-2}{2} E(r)
        +
        r
        \int_{ \S^{n-1}(x,r) }
            |\partial_{r} V|^2 - \theta^2 |\partial_{r}|V||^2
        \,ds
        -
        \frac{r}{2} E'(r).
    \end{align*}
    Solving for $E'(r)$ yields the claim.
    \end{proof}

    \begin{lemma}\label{lemma-alternate-representations-using-the-laplacian}
    For almost every $r$,
    \begin{align*}
        E(r)
        &=
        \frac{k}{2}
        \int_{ \S^{n-1}(x,r) }
            \partial_{r}( |V|^2 )
        \,ds.
    \end{align*}
    \end{lemma}
    \begin{proof}
        The representation of $E(r)$ follows from \eqref{equation-energy-in-terms-of-laplacian}, a smooth approximation of $|V|^2$, and the divergence theorem.
    \end{proof}
    The following spherical-mean identity is well-known for subharmonic functions.
    \begin{lemma}\label{lemma-downstairs-integration-by-parts}
        For almost every $r$, it holds that
        \begin{align*}
            H'(r) = \frac{n-1}{r} H(r) + \int_{ \S^{n-1}(x,r) } \partial_{r}( |V|^2 ) \,ds
        \end{align*}
        and
        \begin{align*}
            H'(r) &= \frac{n-1}{r} H(r) + \int_{ B(x,r) } \Delta(|V|^2) \,dz = \frac{n-1}{r}H(r) + \frac{2}{k} E(r).
        \end{align*}
    \end{lemma}
    \begin{proof}
        It suffices to prove the claim for smooth and nonnegative subharmonic $\varphi$ in place of $|V|^2$ by a standard convolution argument. We have that
        \begin{align*}
            h(r) = \int_{\S^{n-1}(x,r)} \varphi \,ds
            \quad\text{and}\quad
            \delta(r) = \int_{B(x,r)} \Delta\varphi \,dz.
        \end{align*}
        Then, using the radial dilation $T_{x,r} \colon y \mapsto x + r y$,
        \begin{align*}
            h(r) = \int_{\S^{n-1}} (\varphi \circ T_{x,r}) r^{n-1} \,ds
        \end{align*}
        so differentiation with respect to $r$ leads to
        \begin{align*}
            h'(r) 
            &= \frac{n-1}{r} h(r) + \int_{ \S^{n-1} } (\partial_{r}\varphi)( T_{x,r} ) r^{n-1} \,ds
            \\
            &= \frac{n-1}{r} h(r) + \int_{ \S^{n-1}(x,r) } (\partial_{r}\varphi) \,ds
            = \frac{n-1}{r} h(r) + \delta(r),
        \end{align*}
        where the last equality uses the divergence theorem.
    \end{proof}

    We need one last identity before the main monotonicity proof.
    \begin{lemma}\label{lemma-the-final-formula}
    Almost every $r \in \{ e > 0 \}$ satisfies
    \begin{gather*}
        \int_{ \S^{n-1}(x,r) }
            |\partial_{r} V|^2 - \theta^2 |\partial_{r}|V||^2
        \,ds
        -
        \frac{ n }{ 2 }
        \frac{ \mathcal{N}(r) }{ r }
        E(r)
        \\
        =
            \int_{ \S^{n-1}(x,r) }
                \left| \partial_{r}V - \frac{n}{2}\frac{\mathcal{N}(r)}{r} V \right|^2
                -
                \theta^2\left| \partial_{r}|V| - \frac{n}{2}\frac{\mathcal{N}(r)}{r} |V| \right|^2
            \,ds.
    \end{gather*}
    \end{lemma}
    \begin{proof}
    By the chain rule, we conclude
    \begin{align*}
        |V| \partial_{r}( |V| )
        =
        \frac{1}{2}
        \partial_{r}( |V|^2 )
        &=
        \langle V, \partial_{r}V \rangle,
    \end{align*}
    almost everywhere, so
    \begin{align*}
        | \partial_{r}V|^2 - \lambda \frac{\mathcal{N}(r)}{r} \partial_{r}( |V|^2 )
        &=
        \left| \partial_{r}V - \lambda \frac{\mathcal{N}(r)}{r} V \right|^2
        -
        \left( \lambda \frac{\mathcal{N}(r)}{r} \right)^2 |V|^2,
        \\
        | \partial_{r}|V| |^2 - \lambda \frac{ \mathcal{N}(r) }{ r } \partial_{r}( |V|^2) 
        &=
        \left| \partial_{r}|V| - \lambda \frac{ \mathcal{N}(r) }{ r } |V| \right|^2
        -
        \left( \lambda \frac{ \mathcal{N}(r) }{ r} \right)^2 |V|^2
    \end{align*}
    almost everywhere. Hence, recalling Lemma~\ref{lemma-alternate-representations-using-the-laplacian} and the definitions \eqref{equation-dimensional-constant} and \eqref{equation-dimensional-constant-two}, almost every $r \in \{ e > 0 \}$ satisfies
    \begin{gather*}
            \int_{ \S^{n-1}(x,r) }
                |\partial_{r} V|^2 - \theta^2 |\partial_{r}|V||^2
            \,ds
            \\
            =
            \int_{ \S^{n-1}(x,r) }
                \left| \partial_{r}V - \lambda \frac{\mathcal{N}(r)}{r} V \right|^2
                -
               \theta^2\left| \partial_{r}|V| - \lambda \frac{\mathcal{N}(r)}{r} |V| \right|^2
            \,ds
            \\
            +
            \left( 1 - \theta^2 \right)
            \lambda \frac{\mathcal{N}(r)}{r}
            \frac{2}{ k } E(r)
            - 
            \left( 1 - \theta^2 \right) \left( \lambda \frac{\mathcal{N}(r)}{r} \right)^2 H(r)
            \\
            =
            \int_{ \S^{n-1}(x,r) }
                \left| \partial_{r}V - \lambda  \frac{\mathcal{N}(r)}{r} V \right|^2
                -
                \theta^2\left| \partial_{r}|V| - \lambda  \frac{\mathcal{N}(r)}{r} |V| \right|^2
            \,ds
            \\
            +
            \left( 1 - \theta^2 \right)
            \left( \lambda \frac{\mathcal{N}(r)}{r} \right)
            \left(
                \frac{2}{ k } E(r)
                -
                \lambda \frac{\mathcal{N}(r)}{r} H(r)
            \right).
    \end{gather*}
    Recalling that $k = 1-\theta^2$ and letting $\lambda = n/2$, we conclude the claim from \eqref{equation-energy-in-terms-of-laplacian} and \eqref{equation-alternate-expression-of-frequency}.
    \end{proof}

    \subsection{The monotonicity of the frequency}\label{section-the-monotonicity-of-the-frequency}
    We prove Theorem~\ref{theorem-frequency-function-monotonicity} now. To simplify notation, we omit $x$ and $F$ from the computations below.
    
    We recall from Remark~\ref{remark-monotonicity-formula-deconstrunction} that it suffices to consider the monotonicity in the interval $\{e > 0\}$. In this set, by the product rule of absolutely continuous functions,
    \begin{align*}
        \frac{ \mathcal{N}'(r) }{ \mathcal{N}(r) }
        &=
        \frac{1}{r}
        +
        \frac{ E'(r) }{ E(r) }
        -
        \frac{ H'(r) }{ H(r) }.
    \end{align*}
    Next, by Lemmas~\ref{lemma-key-identities-for-the-energy} and \ref{lemma-downstairs-integration-by-parts}, for almost every $r$,
    \begin{align*}
        \frac{ \mathcal{N}'(r) }{ \mathcal{N}(r) }
        &=
        2
        \frac{ 
            \int_{ \S^{n-1}(x,r) }
                |\partial_{r} V|^2 - \theta^2 |\partial_{r}|V||^2
            \,ds
        }{ E(r) }
        -
        \frac{2}{k} \frac{ E(r) }{ H(r) }
        \\
        &=
        \frac{2}{ E(r) }
        \left(
            \int_{ \S^{n-1}(x,r) }
                |\partial_{r} V|^2 - \theta^2 |\partial_{r}|V||^2
            \,ds
            -
            \frac{ n }{ 2 }
            \frac{ \mathcal{N}(r) }{ r }
            E(r)
        \right).
    \end{align*}
    By Lemma~\ref{lemma-the-final-formula}, almost every such $r$ satisfies
    \begin{align*}
        \frac{ \mathcal{N}'(r) }{ \mathcal{N}(r) }
        &=
        \frac{2}{ E(r) }
        \left(
            \int_{ \S^{n-1}(x,r) }
                \left| \partial_{r}V - \frac{n}{2}\frac{\mathcal{N}(r)}{r} V \right|^2
                -
                \theta^2\left| \partial_{r}|V| - \frac{n}{2}\frac{\mathcal{N}(r)}{r} |V| \right|^2
            \,ds
        \right).
    \end{align*}
    We have
    \begin{gather*}
        \left| \partial_{r}|V| - \frac{n}{2}\frac{\mathcal{N}(r)}{r} |V| \right|^2
        \\
        =
        |V|^2 \left( 
            \left\langle \frac{ \partial_r V }{ |V|}, \frac{V}{|V|} \right\rangle^2 
            - 
            n\frac{\mathcal{N}(r)}{r} 
            \left\langle \frac{ \partial_r V }{ |V|}, \frac{V}{|V|} \right\rangle
            +
            \left( \frac{n}{2} \frac{ \mathcal{N}(r) }{ r } \right)^2
        \right)
        \\
        \leq
        |V|^2 \left( 
            \left\langle \frac{ \partial_r V }{ |V|}, \frac{ \partial_r V }{ |V|} \right\rangle^2 
            - 
            n\frac{\mathcal{N}(r)}{r} 
            \left\langle \frac{ \partial_r V }{ |V|}, \frac{V}{|V|} \right\rangle
            +
            \left( \frac{n}{2} \frac{ \mathcal{N}(r) }{ r } \right)^2
        \right)
        \\
        =
        \left| \partial_{r}V - \frac{n}{2}\frac{\mathcal{N}(r)}{r} V \right|^2
    \end{gather*}
    everywhere in $U_F$. In $S_F$, the locality of Sobolev differentiation also implies that
    \begin{align*}
        \left| \partial_{r}|V| - \frac{n}{2}\frac{\mathcal{N}(r)}{r} |V| \right|^2
        =
        0
        =
        \left| \partial_{r}V - \frac{n}{2}\frac{\mathcal{N}(r)}{r} V \right|^2
    \end{align*}
    almost everywhere on $S_F$. 
    
    As $1-\theta^2 = k$, we conclude
    \begin{gather}\label{equation-characterization-of-the-derivative}
        \frac{ 2 }{ E(r) }
            \int_{ \S^{n-1}(x,r) }
                \left| \partial_{r}V - \frac{n}{2}\frac{\mathcal{N}(r)}{r} V \right|^2
            \,ds
        \geq
        \frac{ \mathcal{N}'(r) }{ \mathcal{N}(r) }
        \\ \notag
        \geq
        \frac{ 2k }{ E(r) }
            \int_{ \S^{n-1}(x,r) }
                \left| \partial_{r}V - \frac{n}{2}\frac{\mathcal{N}(r)}{r} V \right|^2
            \,ds
    \end{gather}
    almost everywhere on $\{ e > 0 \}$. The claimed monotonicity follows.

    The following subsections contain applications of the monotonicity, Theorem~\ref{theorem-frequency-function-monotonicity}, and the doubling estimate, Theorem~\ref{theorem-first-fundamental-theorem-of-nevanlinna-theory-and-frequency}.

    \subsection{Proof of Theorem~\ref{theorem-unique-continuation-principle-open-set}}\label{section-unique-continuation-principle}
    The doubling property, Theorem~\ref{theorem-first-fundamental-theorem-of-nevanlinna-theory-and-frequency}, leads to an immediate proof of the unique continuation property of conformal curves.
    \begin{proof}[Proof of Theorem~\ref{theorem-unique-continuation-principle-open-set}]
        If $\mathrm{int}(S_F) = S_F \neq \emptyset$, the equality $S_F = \Omega$ follows, and hence $F$ is constant by the fundamental theorem of calculus. Next, suppose by contradiction that $\emptyset \neq \Omega \cap \overline{U}_F \cap \overline{ \mathrm{int}(S_F) }$. By considering $y \in \Omega \cap \overline{U}_F \cap \overline{ \mathrm{int}(S_F) }$ and $x \in \mathrm{int}(S_F)$ close enough to $y$, we lose no generality in supposing
        \begin{align*}
            r_0 \coloneqq \mathrm{dist}( x, U_F ) < \mathrm{dist}( x, \R^n \setminus \Omega ) = r(x).
        \end{align*}
        By Remark~\ref{remark-monotonicity-formula-deconstrunction}, it holds that
        \begin{align*}
            h_{F}(x, r_0) = 0
            \quad\text{and}\quad
            h_{F}(x,r) > 0
            \quad\text{for $r_0 <r<r(x)$.}
        \end{align*}
        We let $r_0 < r_1 < r(x)$ for which
        \begin{align*}
            n
            \mathcal{N}_{F}(x,r_1) \log\left( \frac{r_1}{r_0} \right)
            \leq
            \frac{1}{2}.
        \end{align*}
        Then Theorem~\ref{theorem-first-fundamental-theorem-of-nevanlinna-theory-and-frequency} implies
        \begin{align*}
            h_{F}(x,r)
            \leq
            \frac{1}{2}
            h_{F}(x,r)
        \end{align*}
        for every $r_0 < r < r_1$. Hence, $h_{F}(x,r) = 0$, leading to a contradiction.
    \end{proof}
    This leads to the following convenient lemma.
    \begin{lemma}\label{lemma-consequence-of-monotonicity}
        Let $F \colon \R^n \supset \Omega \to ( \R^m, \omega )$ be a conformal curve for $3 \leq n \leq m$. If $F$ is non-constant, then $h_{F}(x,r) > 0$ for every $x \in \Omega$ and $0 < r < r(x)$.
    \end{lemma}
    \begin{proof}
        If $h_{F}(x,r) = 0$, then $V_F = 0$ on $\S^{n-1}(x,r)$. By the maximum principle of the subharmonic function $|V_F|^2$, we conclude that $DF = 0$ on $B(x,r)$. Hence $F|_{B(x,r)}$ is constant. By the unique continuation principle, Theorem~\ref{theorem-unique-continuation-principle-open-set}, $F$ itself is constant.
    \end{proof}

    \subsection{The doubling property}\label{section-doubling-property}
    Theorem~\ref{theorem-first-fundamental-theorem-of-nevanlinna-theory-and-frequency} leads to the following quantitative doubling property of the energy.
    \begin{corollary}\label{corollary-doubling-property-of-the-energy}
        Let $\mathcal{N}_{0}>0$ be such that
        \begin{align*}
            \mathcal{N}_{F}(x,r_0) < \mathcal{N}_{0}
        \end{align*}
        and $\lambda_0 > 1$ such that
        \begin{align*}
            c_{\lambda_0}
            =
            n
            \mathcal{N}_{0}
            \log\left( \lambda_0 \right)
            <
            1.
        \end{align*}
        Then it holds that
        \begin{align*}
            \int_{ B(x,\lambda_{0} r) } |V_F|^2 \,dz 
            \leq
            \frac{ \lambda_{0} }{ 1 - c_{\lambda_{0}} }
            \int_{ B(x, r) } |V_F|^2 \,dz
        \end{align*}
        for every $0 < r < r_0$.
    \end{corollary}
    \begin{proof}
        By Fubini's theorem and Theorem~\ref{theorem-first-fundamental-theorem-of-nevanlinna-theory-and-frequency},
        \begin{align*}
            \int_{ B(x, \lambda_{0} r) } |V_F|^2 \,dz
            &=
            \int_{ 0 }^{ \lambda_{0} r } H_{F}(x,s) \,ds
            =
            \lambda_{0}
            \int_{ 0 }^{ r } H_{F}(x,\lambda_{0}s) \,d s
            \\
            &\leq
            \frac{ \lambda_{0} }{ 1 - c_{\lambda_{0}} }
            \int_{0}^{r}
                H_{F}(x,s)
            \,ds
            =
            \frac{ \lambda_{0} }{ 1 - c_{\lambda_{0}} }
            \int_{ B(x, r) } |V_F|^2 \,dz.
        \end{align*}
    \end{proof}
    Let $\ell = \ell(\lambda_{0}) \in \N$ be the least integer for which $\lambda_{0}^{\ell} \geq 2$. We let
    \begin{equation*}
        C_0
        =
        \left(
            \frac{ \lambda_{0} }{ 1 - c_{\lambda_{0}} }
        \right)^{ \frac{\ell}{n} },
    \end{equation*}
    noting the doubling property
    \begin{align}\label{equation-doubling-property-of-the-energy}
        \int_{ B(x, 2r) } |V_F|^2 \,dz
        \leq
        C_{0}^n
        \int_{ B(x,r) } |V_F|^2 \,dz.
    \end{align}

\subsection{A modulus of continuity}\label{section-modulus-of-continuity}
    Using \eqref{equation-doubling-property-of-the-energy}, we obtain a uniform Lipschitz estimate for $F$ itself.
    \begin{proposition}\label{proposition-energy-controls-diameter}
        For every $0 < 2r < r_0$, it holds that
        \begin{align*}
            |F(y)-F(z)|
            \leq
            C_0
            c_n
            \left(
                \aint{ B(x,r) } |V_F|^2 \,dz
            \right)^{ \frac{1}{n} }
            |y-z|
        \end{align*}
        for every $y,z \in B( x, r )$. Here $c_n$ depends only on $n$.
    \end{proposition}
    \begin{proof}
        By \cite{ikonen-2026-pushforward-of-currents-under-sobolev-maps} (see \cite[Proposition~7.1]{ikonen-pankka-2024-liouvilles-theorem-in-calibrated-geometries}), it holds that
        \begin{align*}
            |F(y)-F(z)|
            \leq
            c_n'
            \left(
                \aint{ B(x,2r) } |V_F|^2 \,dz
            \right)^{ \frac{1}{n} }
            |y-z|
        \end{align*}
        for every $y, z \in B( x, r )$. Since
        \begin{align*}
            \left(
                \aint{ B(x,2r) } |V_F|^2 \,dz
            \right)^{ \frac{1}{n} }
            \leq
            \frac{ C_0 }{ 2 }
            \left(
                \aint{ B(x,r) } |V_F|^2 \,dz
            \right)^{ \frac{1}{n} },
        \end{align*}
        the claim follows by setting $c_{n} = c_{n}' / 2$.
    \end{proof}
    
\subsection{First-order estimates}\label{section-first-order-estimates}
    The regularity theory of $n$-harmonic maps also leads to the following uniform $C^{1,\alpha}_{loc}$-estimates using Corollary~\ref{corollary-doubling-property-of-the-energy}.
    \begin{proposition}\label{proposition-c-1-alpha-estimates}
        There exists a universal $\alpha = \alpha(n) \in (0,1)$ with the following property. For every $0 < 8r < r_0$, it holds that
        \begin{align*}
            \sup_{ \substack{y \neq z \\ y,z \in B(x,r)} }
                \frac{ |DF(y)-DF(z)| }{ |y-z|^{\alpha} }
            &\leq
            \frac{ (2C_0^3) c_n C_{n} }{ r^{\alpha} }
            \left(
               h_{F}(x,r)
            \right)^{ \frac{1}{n} }.
        \end{align*}
        Here $C_{n}$ depends only on $n$.
    \end{proposition}
    \begin{proof}
        The standard $C^{1,\alpha}$-regularity estimate by Uhlenbeck \cite[Theorems~(3.2) and (5.4)]{uhlenbeck-1977-regularity-for-a-class-of-non-linear-elliptic-systems} together with Proposition~\ref{proposition-energy-controls-diameter} and \eqref{equation-doubling-property-of-the-energy} give that
        \begin{align*}
            \sup_{ \substack{y \neq z \\ y,z \in B(x,r)} }
                \frac{ |DF(y)-DF(z)| }{ |y-z|^{\alpha} }
            &\leq
            \frac{ C_{n}' }{ r^{1+\alpha} }
            \left( \aint{B(x,4r)} |DF|^n \,dz \right)^{\frac{1}{n}}
            \\
            &\leq
            \frac{ C_{n}' }{ r^{1+\alpha} }
            \frac{ C_0^{3} }{ 2^{2} }
            c_n
            \left(
               \aint{ B(x,r) } |V_F|^2 \,dz
            \right)^{ \frac{1}{n} }.
        \end{align*}
        By the isoperimetric inequality for conformal curves, cf.~\cite[Theorem~4.1]{ikonen-pankka-2024-liouvilles-theorem-in-calibrated-geometries}, and Hölder's inequality, we have
        \begin{align*}
            \left(
                \aint{ B(x,r) } |V_F|^2 \,dz
            \right)^{ \frac{1}{n} }
            &\leq
            \left(
               \aint{ \S^{n-1}(x,r) } |V_F|^2 \,ds
            \right)^{ \frac{1}{n} }.
        \end{align*}
        Combining the estimates leads to the assertion.
    \end{proof}
    \subsection{The basic properties of the frequency}\label{section-basic-properties-of-the-frequency}
    Let $F \colon \R^n \supset \Omega \to ( \R^m, \omega )$ be a conformal curve for a constant-coefficient calibration $\omega$ and $3 \leq n \leq m$.

    If $F$ is non-constant, Lemma~\ref{lemma-consequence-of-monotonicity} shows that $h_{F}(x,r) > 0$ for every $x \in \Omega$ and $0 < r < r(x)$. Using this fact, we note the following elementary result.
    \begin{lemma}\label{lemma-continuity-of-frequency}
        The frequency $(x,r) \mapsto \mathcal{N}_{F}(x,r)$ is continuous in the open set $$D = \{ (y,s) \in \Omega \times (0,\infty) \colon r(y) > s \}.$$
    \end{lemma}
    \begin{proof}
        If $F$ is non-constant, then Lemma~\ref{lemma-consequence-of-monotonicity} implies $D = \{ (y,s) \in D \colon h_{F}(y,s) > 0 \}$. By the $C^{1}$-regularity of $F$, the function $(y,s) \mapsto h_{F}(y,s)$ is continuous in $D$. Since also $(y,s) \mapsto e_{F}( y, s )$ is continuous by dominated convergence and the local integrability of $\Delta( |V_F|^2 )$, the claimed continuity follows.

        On the other hand, if $F$ is constant, clearly the claim holds as well.
    \end{proof}

    \section{Radially homogeneous conformal curves}\label{section-radially-homogeneous-conformal-curves}
    For this section, we fix a conformal curve $F \colon \R^n \supset \Omega \to ( \R^m, \omega )$ for a calibration $\omega \in \Lambda^{n} \R^m$ and $3 \leq n \leq m$. We analyze conformal curves for which $r \mapsto \mathcal{N}_{F}(x,r)$ is constant in an interval.

    \subsection{Recognizing radial homogeneity}\label{section-recognizing-radial-homogeneity}
    In this section, we prove Theorem~\ref{theorem-introduction-the-main-frequency-structure-result}. In fact, we prove a slightly stronger version of the statement. To motivate the precise claim, consider a conformal curve satisfying $F(x+rv) = F(x) + r^{\gamma+1} h(v)$ whenever $x+rv \in \Omega$, $r >0$, and $|v| = 1$, for some function $h \colon \S^{n-1} \to \R^m$, $\gamma \geq 0$, and $x \in \Omega$. Then a direct computation shows $\mathcal{N}_{F}(x,r) = \gamma \geq 0$ for every $0 < r < r(x)$. The following theorem can be understood as a converse of this result.
    \begin{theorem}\label{theorem-the-main-frequency-structure-result}
        Let $x \in \Omega$ and $0 \leq s_0 < r_0 \leq r(x)$. If $\mathcal{N}_{F}(x,r) = \mathcal{N}_0$ for every $s_0 < r < r_0$, then there exists a real-analytic $h \colon \S^{n-1} \to \R^{m}$ such that
        $$F(x+rv) = F(x) + r^{ \mathcal{N}_{0} + 1 } h(v) \quad\text{whenever $x+rv \in \Omega$.}$$ Furthermore,
        \begin{enumerate}
            \item $|h| = C$ for a constant $C \geq 0$;
            \item $F^{*}g_{euc} = C^2 \left( \mathcal{N}_0 + 1 \right)^2 r^{ 2\mathcal{N}_0 } g_{euc}$ and $h^{*}g_{euc} = C^2 \left( \mathcal{N}_{0}+1\right)^2 g_{\S^{n-1}}$;
            \item $\Delta_{ \S^{n-1} } h = -(n-1) \left( \mathcal{N}_{0} + 1 \right)^2 h$.
        \end{enumerate}
        In particular, the frequency satisfies $\mathcal{N}_{F}(x,r) = \mathcal{N}_0$ for every $0 < r < r(x)$.
    \end{theorem}
    \begin{remark}
        The proof shows that if $\mathcal{N}_{F}(x,r) = 0$ for some $0 < r < r(x)$, then $F|_{\Omega}$ is affine. This is a stronger conclusion than stated in Theorem~\ref{theorem-the-main-frequency-structure-result}. Theorem~\ref{theorem-the-main-frequency-structure-result} is significantly stronger than Theorem~\ref{theorem-introduction-the-main-frequency-structure-result}. Indeed, the frequency being constant in an interval is strong enough to describe the whole map $F \colon \Omega \to ( \R^m, \omega )$.
    \end{remark}
    \begin{remark}\label{remark-spectrum-of-the-laplacian}
        Theorem~\ref{theorem-the-main-frequency-structure-result} is closely related to the classification of isometric immersions $\S^{n-1}(0,c) \rightarrow \S^{m-1}$ due to Takahashi \cite{takahashi-1966-minimal-immersions-of-riemannian-manifolds} and do Carmo--Wallach \cite{do-carmo-wallach-1971-minimal-immersions-of-spheres-into-spheres}; see \cite[Theorem~3]{takahashi-1966-minimal-immersions-of-riemannian-manifolds} in particular. Indeed, using the notation of Theorem~\ref{theorem-the-main-frequency-structure-result}, the map $$\S^{n-1} \to \S^{m-1}\left( 0, \left( \frac{1}{ \mathcal{N}_{0}+1 } \right) \right), \quad v \mapsto \frac{ h(v) }{  C \left( \mathcal{N}_0 + 1 \right) }$$ is an isometric immersion. Equivalently, $$g \colon \S^{n-1}\left( 0, \mathcal{N}_0 + 1 \right) \rightarrow \S^{m-1}( 0, 1 ), \quad v \mapsto \frac{ h\left( \frac{v}{ \mathcal{N}_{0} + 1 } \right) }{ C } $$ defines an isometric immersion.
        
        The spectrum of the spherical Laplacian (see e.g. \cite[Lemma~22, page~39]{muller-1966-spherical-harmonics}) gives that either $F$ is constant on $\Omega$ or
        \begin{align*}
            (n-1) \left( \mathcal{N}_{0} + 1 \right)^2 = \ell(\ell+n-2)
        \end{align*}
        for $\ell \in \{1,2,\dots\}$. We recall that $\ell$ is the degree of radially homogeneous harmonic polynomial $P \colon \R^n \to \R^m$ satisfying $P|_{ \S^{n-1} } = h$. Thus, in particular,
        \begin{align*}
            \mathcal{N}_{0} + 1
            &=
            \left( \frac{ \ell(\ell+n-2) }{ n - 1 } \right)^{ \frac{1}{2} }.
        \end{align*}
        Hence, the possible degrees of non-constant $F$ are
        \begin{align*}
            1, \left( \frac{2n}{n-1} \right)^{ \frac{1}{2} }, \left( \frac{3(n+1)}{n-1} \right)^{ \frac{1}{2} }, \dots.
        \end{align*}
        In fact, it is clear that $F$ is affine if and only if $\ell=1$, or, equivalently, if $rv \mapsto F(x+rv)= F(x) + P(rv)$ for $x+rv \in \Omega$.
        
        Moreover, a simple computation shows that
        \begin{align}\label{equation-normalized-harmonics}
            h_{F}(x,1)
            =
            n^{ \frac{n}{2}{} }
            \Rightarrow
            C 
            &= 
            \frac{ 1 }{ \mathcal{N}_{0} + 1 }
            =
            \frac{ 1 }{ \left( \frac{ \ell(\ell+n-2) }{ n - 1 } \right)^{ \frac{1}{2} } }.
        \end{align}
        We see from Theorem~\ref{theorem-the-main-frequency-structure-result} (2) that \eqref{equation-normalized-harmonics} is equivalent to $h$ being an isometric immersion. Equivalently, the restriction $F|_{ \overline{B} }$ is $1$-Lipschitz.
    \end{remark}

    In order to simplify our notation for the proof, we omit the dependence on $x$ and $F$ in our notations. To prepare for the proof of Theorem~\ref{theorem-the-main-frequency-structure-result}, we first show the following proposition.
    \begin{proposition}\label{proposition-recognizing-constant-frequency}
        It holds that $\mathcal{N}(r) = 0$ for some $0<r<r(x)$ if and only if $F$ is affine on $\Omega$.
        
        Let $\mathcal{N}_0 \in (0,\infty)$ and $0 \leq s_0 < r_0 \leq r(x)$. It holds that $\mathcal{N}(r) = \mathcal{N}_0$ for every $s_0<r<r_0$ if and only if
        \begin{align}\label{equation-radial-homogeneity-of-V}
            \partial_{r}V = \frac{n}{2} \frac{ \mathcal{N}(r) }{ r }V
        \end{align}
        almost everywhere in $\S^{n-1}( x, r )$ for almost every $s_0 < r < r_0$.
    \end{proposition}
    \begin{proof}
        First, if $F$ is constant, then clearly $\mathcal{N}(r) = 0$ for every $0 < r < r(x)$. The claim is clear in this case.

        We assume for the rest of the claim that $F$ is non-constant. By Lemma~\ref{lemma-consequence-of-monotonicity}, we have $\{ h > 0 \} = (0,r(x))$.

        We consider two cases. In the first case, we suppose that $e(r) = 0$ for some $r \in (0, r(x))$. By definition, $\Delta( |V|^2 ) = 0$ almost everywhere on $B(x,r)$. By Corollary~\ref{corollary-representation-of-the-laplacian}, it holds that $|DV| = 0$ almost everywhere on $B(x,r)$, i.e., that $V$ is a constant on $B(x,r)$. Since $A \mapsto |A|^{ \frac{n-2}{2} } A$ is a homeomorphism on $\R^{ m \times n }$, it follows that $DF$ is constant on $B(x,r)$. In particular, there exists an affine map $A \colon \R^n \to \R^m$ for which $A|_{B(x,r)} = F|_{ B(x,r) }$. 
        
        We claim that $A|_{\Omega} = F$. By construction, $A$ is not a constant map. Moreover, the relatively closed set $C = \{ y \in \Omega \colon DA(y) = DF(y) \}$ contains $B(x,r)$. Since $DA$ has full rank, it follows that $C \subset U_F$. By the unique continuation principle for real-analytic functions, $F|_{W} = A|_{W}$ in the connected component $W$ of $U_F$ containing $B(x,r)$. By openness of $W$, it clearly follows that $W \subset C$. If $y \in \Omega$ is an accumulation point of $W$, then $y \in C$. Since $C \subset U_F$ and $U_F$ is open, we conclude $B(y, \delta ) \subset U_F$ for small enough $\delta > 0$. By definition, $B(y,\delta) \subset W$. We conclude that $W$ is relatively closed in $\Omega$. By connectivity, $W = \Omega$ follows. This proves the claim. Now, conversely, if $F = A|_{\Omega}$ for a non-constant affine map, clearly $\mathcal{N}(r) = 0$ for every $0 < r < r(x)$.

        In the final case, we have $\{ e > 0 \} = (0,r(x))$. We claim that $\mathcal{N}(r) = \mathcal{N}_0$ for every $r \in (s_0,r_0)$ if and only if \eqref{equation-radial-homogeneity-of-V} holds. However, this is immediate from \eqref{equation-characterization-of-the-derivative}.
    \end{proof}

    \begin{lemma}\label{lemma-auxiliary-radial-homogeneity}
        Let $\mathcal{N}_0 \in [0,\infty)$. It holds that $\mathcal{N}(r) = \mathcal{N}_{0}$ for every $s_0<r<r_0$ if and only if $DF(x+rv) = \left( \frac{r}{s} \right)^{  \mathcal{N}_{0} } DF(x+sv)$ for every $s_0 < s, r < r_0$ and $v \in \S^{n-1}$.
    \end{lemma}
    \begin{proof}
        The case when $\mathcal{N}_0 = 0$ is clear from Proposition~\ref{proposition-recognizing-constant-frequency}. It remains to consider the case $\{ e > 0 \}= (0,r(x))$. Since $A \mapsto |A|^{ \frac{n-2}{2} } A$ is invertible in $\R^{m \times n}$, the conclusion for $DF$ in the claim is equivalent to
        \begin{align}\label{equation-homogeneity}
            V\left( x+r v \right) =\left( \frac{r}{s} \right)^{\frac{n}{2}\mathcal{N}_{0}} V( x+s v )
            \quad\text{for $r,s \in (s_0,r_0)$.}
        \end{align}
        When \eqref{equation-homogeneity} holds, the radial differentiation stated in \eqref{equation-radial-homogeneity-of-V} clearly holds, and thus $\mathcal{N}(r) = \mathcal{N}_0$ holds. Conversely, we first note that $V$ is locally Hölder continuous and in $W^{1,2}_{loc}$, so the restriction of $V$ to almost every radial ray is in $W^{1,2}_{loc}( (0,r_0), \R^{ m \times n } )$. Now, for almost any such direction $v \in \mathbb{S}^{n-1}$, the continuous matrix field $a_v(r) = V(x+rv)$ satisfies
        \begin{align*}
            a_{v}'(r) = \frac{\frac{n}{2}\mathcal{N}_{0}}{r} a_{v}(r)
        \end{align*}
        for almost every $r \in (s_0,r_0)$ by Proposition~\ref{proposition-recognizing-constant-frequency}, leading to  
        \begin{align*}
            a_{v}(r) = \left( \frac{r}{s} \right)^{\frac{n}{2}\mathcal{N}_{0}} a_{v}( s )
        \end{align*}
        for every $s,r \in (s_0,r_0)$. By continuity of $V$ and density of such directions, the claim follows.
    \end{proof}
    \begin{proof}[Proof of Theorem~\ref{theorem-the-main-frequency-structure-result}]
        Again, if $\mathcal{N}_{F}(x,r) = 0$ for some $0 < r < r(x)$, we have that $F$ is affine by Proposition~\ref{proposition-recognizing-constant-frequency}. Hence, we may suppose that $F$ is not affine on any open subset of $\Omega$. In particular, $\{ e > 0 \} = (0,r(x))$. Now, we may fix $s_0 < s < r_0$ and denote
        \begin{align*}
            H(v)[w]
            =
            \frac{ 1 }{ \mathcal{N}_0 + 1 }
            \left(
                \frac{ DF(x+sv)[w] }{ s^{ \mathcal{N}_0} }  
            \right)
            \quad\text{and}\quad
            h(v)
            =
            H(v)[v]
        \end{align*}
        for $v \in \S^{n-1}$ and $w \in \R^{n}$. Using the notation and conclusion of Lemma~\ref{lemma-auxiliary-radial-homogeneity} and $C^1$-regularity of $F$, we have
        \begin{align*}
            DF(x+rv)
            =
            \left( \mathcal{N}_0 + 1 \right)
            r^{ \mathcal{N}_{0} }
            H(v).
        \end{align*}
        Clearly this equality holds for $r = s_0$ as well. Thus, integrating along a radial line segment from $x+s_0v$ to $x+rv$ leads to
        \begin{align}\label{equation-homogeneous-representation-of-F}
            F(x+rv)
            -
            F(x+s_0v)
            =
            \left( r^{ \mathcal{N}_{0} + 1 } - s_0^{ \mathcal{N}_{0}+1 } \right)
            h(v).
        \end{align}
        Now, let $\theta(t) = \cos(t)v+\sin(t)w$ for orthogonal and unit $v$ and $w$. Differentiating the representation \eqref{equation-homogeneous-representation-of-F} along $\gamma(t) = x+r \theta(t)$ at $t = 0$ implies
        \begin{align}\label{equation-reformulation-of-the-matrix-field}
            (Dh)(v)[w] = \left( \mathcal{N}_{0}+1\right)H(v)[w].
        \end{align}
        Denoting by $T$ the usual matrix transpose, the conformality of $F$ yields
        \begin{align}\label{equation-rephrasing-of-the-conformality-condition}
            H^T H
            =
            \frac{ |H|^2 }{ n } I_{ n \times n }
        \end{align}
        hence
        \begin{align}\label{equation-a-radial-consequence-of-conformality}
            \langle (DH)(v)[w], h(v) \rangle = 0.
        \end{align}
        Combining \eqref{equation-reformulation-of-the-matrix-field} and \eqref{equation-a-radial-consequence-of-conformality} leads to the fact that $|H(v)[v]|^2/n = |h(v)|^2 = C^2$ for some $C \geq 0$ for every $v \in \S^{n-1}$ as claimed. Hence, by \eqref{equation-rephrasing-of-the-conformality-condition},
        \begin{align}\label{equation-rephrasing-of-the-conformality-condition-improved}
            H^T H
            =
            C^2 I_{ n \times n }.
        \end{align}
        Here \eqref{equation-reformulation-of-the-matrix-field} and \eqref{equation-rephrasing-of-the-conformality-condition-improved} imply that
        \begin{align}\label{equation-conformality-immersion-property-of-h}
            (Dh)^T (Dh)
            =
            \frac{ |Dh|^{2} }{ n-1 } I_{ T \S^{n-1} }
            =
            C^2 \left( \mathcal{N}_0 + 1 \right)^2 I_{ T \S^{n-1} }.
        \end{align}
        Here \eqref{equation-conformality-immersion-property-of-h} implies the pullback property for $h$.
        
        Lemma~\ref{lemma-consequence-of-monotonicity} implies $C>0$, so $|H(v)[v]| = |h(v)| = C > 0$. In particular, the differential of $F$ does not vanish on the annulus $A = B( x, r_0 ) \setminus \overline{B}( x, s_0 )$, i.e., $A \subset U_F$. As $F$ is real-analytic on $U_F$, we conclude that $h$ is real-analytic. We define the map
        \begin{align*}
            G(x+rv) = F(x) + r^{ \mathcal{N}_0+1 } h(v)
            \quad\text{whenever $x+rv \in \Omega$.}
        \end{align*}
        It holds that $G$ is real-analytic on $\Omega \setminus \{x\}$ and $A \subset C = \{ DG = DF \}$.
        
        In fact, we claim $F = G$. We prove this by showing that $C$ is open and relatively closed in $\Omega$. First, since $DG$ has full rank in $\Omega \setminus \{x\}$, we have $C \setminus \{x\} \subset U_F \setminus \{x\}$. Also, the unique continuation principle for real-analytic functions guarantees that $DG = DF$ in the connected component $W$ of $U_F \setminus \{x\}$ containing $A$. If $y \in \Omega \setminus \{x\}$ is an accumulation point of $W$, then $y \in C$ by construction. Thus $y \in U_F$ as above, implying that $B(y,\delta) \subset U_F \setminus \{x\}$ for small enough $\delta > 0$. Therefore, $B(y,\delta) \subset W$. It follows that $W$ is open and relatively closed in $\Omega \setminus \{x\}$, implying $W = \Omega \setminus \{x\}$ and $C \setminus \{x\} = \Omega \setminus \{x\}$. By continuity, it also holds that $C = \Omega$. Since the differentials of $F$ and $G$ agree on $\Omega$ and $F(x) = G(x)$, we conclude that $F = G$ in $\Omega$.

        We have now verified that $F(x+rv) = F(x) + r^{ \mathcal{N}_{0} + 1 }h(v)$ on $\Omega$. It, therefore, follows that
        \begin{align*}
            ( rv \mapsto F(x+rv) )^{*}g_{euc}
            =
            C^2
            \left(
                \mathcal{N}_0 + 1 
            \right)^2
            r^{ 2\mathcal{N}_0 }
            g_{euc}.
        \end{align*}
        Only the Laplacian claim remains to be verified. To this end, recall from \eqref{equation-conformality-and-n-harmonicity}, for the components $F = \sum_{ k } F_k e_k$, that
        \begin{align*}
            \Delta F_k &= -(n-2) \langle Df, DF_k \rangle,
        \end{align*}
        where $f = \log|DF|$ and $k = 1,2,\dots,m$. But here
        \begin{align*}
            f(x+rv) = \log\left( \sqrt{n} C\left( \mathcal{N}_{0}+1\right) \right) +  \mathcal{N}_{0} \log(r).
        \end{align*}
        Therefore,
        \begin{align*}
            Df(x+rv) 
            =
            \mathcal{N}_{0}
            \frac{ v }{ r }.
        \end{align*}
        Moreover, 
        \begin{align*}
            \langle Df, DF_k \rangle(x+rv)
            &=
            \mathcal{N}_{0}
            \left(\mathcal{N}_{0}+1\right)r^{ \mathcal{N}_{0}-1 }
            \langle e_{k}, h(v) \rangle.
        \end{align*}
        Hence,
        \begin{align*}
            \Delta F(x+rv)
            &=
            -(n-2)
            \mathcal{N}_{0}
            ( \mathcal{N}_{0}+1 )r^{ \mathcal{N}_{0} - 1 }
            h(v)
        \end{align*}
        using the following consequences of radial homogeneity of $DF$:
        \begin{align*}
            (\partial_{r}F)(x+rv)
            &=
            ( \mathcal{N}_{0} + 1 )r^{ \mathcal{N}_{0}  } h(v),
            \quad\text{and}
            \\
            (\partial_{r}\partial_{r}F)(x+rv)
            &=
            ( \mathcal{N}_{0} + 1 )
            \mathcal{N}_{0} r^{ \mathcal{N}_{0} -1 } h(v).
        \end{align*}
        On the other hand, expressing the Euclidean Laplacian in spherical coordinates leads to
        \begin{align*}
            \Delta F(x+rv)
            &=
            (\partial_{r}\partial_{r}F)(x+rv)
            +
            \frac{n-1}{r} (\partial_{r}F)(x+rv)
            +
            \frac{1}{r^2} (\Delta_{ \S^{n-1} } F)(x+rv)
            \\
            &=
            \left(
                \mathcal{N}_0
                +
                (n-1)
            \right)
            ( \mathcal{N}_{0} + 1 )
            r^{ \mathcal{N}_{0} -1 } h(v)
            +
            r^{ \mathcal{N}_{0} - 1 } (\Delta_{ \S^{n-1} } h)(v),
        \end{align*}
        whence
        \begin{align*}
            (\Delta_{ \S^{n-1} } h)(v)
            &=
            -
            (n-1)
            \left( \mathcal{N}_{0} + 1 \right)^{2}
            h(v).
        \end{align*}
    \end{proof}

    \subsection{Proof of Theorem~\ref{theorem-normalized-blow-up}}\label{section-blow-ups-of-curves}
    In this subsection, we consider a non-constant conformal curve $F \colon \R^n \supset \Omega \to ( \R^m, \omega )$ and let $x \in \Omega$. 

    For the proof of Theorem~\ref{theorem-normalized-blow-up}, we need the following rescaling lemma of the frequency.
    \begin{lemma}\label{lemma-rescaling-the-radius}
    For every $\overline{B}(y,s) \subset \{ z \in \R^n \colon x+rz \in \Omega \}$,
    \begin{align*}
        \mathcal{N}_{F_r}(y,s) 
        =
        \mathcal{N}_{F}(x+ry,rs)
        \quad\text{and}\quad
        h_{ F_{r} }(y,s)
        =
        \left( \frac{r}{c_r(x)} \right)^n
        h_{F}(x+ry,rs).
    \end{align*}
    In particular, if $(y,s) = (0,1)$, then
    \begin{align*}
        \mathcal{N}_{F_r}(0,1)
        &=
        \mathcal{N}_{F}(x,r)
        \quad\text{and}\quad
        \frac{ h_{F_r}(0,1) }{ n^{ \frac{n}{2} } }
        =
        1.
    \end{align*}
    \end{lemma}
    \begin{proof}
    We compute that
    \begin{align*}
        ( DF_r )(z) &= \left( \frac{r}{c_r(x)} \right) ( DF )( x + rz )
    \end{align*}
    so
    \begin{align*}
        V_{F_r}(z) &= \left( \frac{r}{c_r(x)} \right)^{ \frac{n}{2} } V_{F}( x + rz  ).
    \end{align*}
    This gives
    \begin{align*}
        \Delta( |V_{F_r}|^2 )(z)
        &=
        \left( \frac{r}{c_r(x)} \right)^{ n }
        \Delta( z \mapsto |V_{F}( x + rz )|^2 )(z)
        \\
        &=
        \left( \frac{r}{c_r(x)} \right)^{ n } r^{2} \Delta( |V_{F}|^2 )( x + rz ),
    \end{align*}
    leading to
    \begin{align*}
        \int_{ B(y,s) } \Delta( |V_{F_r}|^2 ) \,dz
        &=
        \left( \frac{1}{c_{r}(x)} \right)^2
        r^{2}
        \int_{ B(x+ry,rs) } \Delta( |V_{F}|^2 ) \,dz.
    \end{align*}
    Therefore, multiplying by $s/n$ and dividing by $\omega_n s^n$ leads to the formula
    \begin{align*}
        s e_{F_r}(y,s)
        &=
        \left( \frac{r}{c_{r}(x)} \right)^{n}
        (rs) e_{F}(x+ry,rs).
    \end{align*}
    Similarly,
    \begin{align*}
        H_{F_r}(y,s)
        &=
        \int_{ \S^{n-1}(y,s) } |V_{F_r}|^2 \,ds
        =
        \left( \frac{1}{c_r(x)} \right)^{ n }
        r
        H_{F}(x+ry,rs).
    \end{align*}
    Dividing both sides by $n\omega_n s^{n-1}$ leads to
    \begin{align*}
        h_{F_r}(y,s) &= \left( \frac{r}{c_r(x)} \right)^{ n } h_{F}(x+ry,rs)
    \end{align*}
    as claimed. The claimed formula for the frequency follows as well.
    \end{proof}
    \begin{proof}[Proof of Theorem~\ref{theorem-normalized-blow-up}]
        We let $\mathcal{N}_0 > \mathcal{N}_{F}(x)$ and let $\lambda_{0}>1$ such that
        \begin{align*}
            c_{\lambda_0}
            =
            n
            \mathcal{N}_{0}
            \log\left( \lambda_0 \right)
            <
            1.
        \end{align*}
        Let $r_0 >0$ with $\mathcal{N}_{F}(x,r_0) < \mathcal{N}_0$.

        For each $\ell \in \N$, Lemma~\ref{lemma-rescaling-the-radius} implies that
        \begin{align*}
            \mathcal{N}_{F_{r_j}}( 0, \lambda_0^{\ell} )
            =
            \mathcal{N}_{F}( x, r_j \lambda_{0}^{\ell} )
            <
            \mathcal{N}_0
        \end{align*}
        whenever $r_j < r_0 \lambda_{0}^{-\ell}$. Therefore, iterating Theorem~\ref{theorem-first-fundamental-theorem-of-nevanlinna-theory-and-frequency} implies
        \begin{align}\label{equation-local-boundedness-of-boundary-energy}
            h_{F_{r_j}}(0,\lambda_0^{\ell})
            \leq
            \left( \frac{1}{1-c_{\lambda_0}} \right)^{\ell}
        \end{align}
        whenever $r_j < r_0 \lambda_{0}^{-\ell}$.
        
        Proposition~\ref{proposition-energy-controls-diameter} and \eqref{equation-local-boundedness-of-boundary-energy} imply that the sequence $( F_{r_j} )$ has a uniform Lipschitz estimate on compact subsets of $\R^n$. Similarly, the normalization in Proposition~\ref{proposition-c-1-alpha-estimates} and \eqref{equation-local-boundedness-of-boundary-energy} imply boundedness of $( F_{r_j} )$ in $C^{1,\alpha}$-norm over $B(0,\lambda_0^{\ell})$, for each $\ell \in \N$. By passing to a diagonal subsequence and relabelling, we obtain convergence $F_{r_j} \rightarrow F_{\infty}$ uniformly and in $C^1$ on every compact subset of $\R^n$.

        We prove that the frequency of $F_\infty$ is constant at $y = 0$. To this end, we first notice that
        \begin{align*}
            h_{F_\infty}(0,r) = \lim_{ j \to \infty } h_{F_{r_j}}(0,r)
        \end{align*}
        for every $r$ by the $C^{1}$-convergence. In particular, for $r = 1$, we have
        \begin{align}\label{equation-blow-up-analysis-normalization}
            h_{F_\infty}(0,1) = n^{ \frac{n}{2} }.
        \end{align}
        Similarly, if $0 < s < r$, \eqref{equation-first-fundamental-theorem-of-nevanlinna-theory} implies
        \begin{align*}
            h_{F_\infty}(0,r)
            -
            h_{F_\infty}(0,s)
            &=
            \lim_{j\to\infty}
            h_{F_{r_j}}(0,r)
            -
            h_{F_{r_j}}(0,s)
            \\
            &=
            \lim_{j\to\infty}
            n
            \int_{s}^{r}
                \frac{ 1 }{ a }
                \mathcal{N}_{F_{r_j}}(0,a)
                h_{ F_{r_j} }(0,a)
            \,da
            \\
            &=
            n
            \int_{s}^{r}
                \frac{ 1 }{ a }
                \mathcal{N}_{F}(x)
                h_{ F_{\infty} }(0,a)
            \,da.
        \end{align*}
        Indeed,
        \begin{align*}
            \mathcal{N}_{F_{r_j}}(0,a) = \mathcal{N}_{F}(x,r_j a) = \mathcal{N}_{F}(x) + o(1)
        \end{align*}
        and
        \begin{align*}
            \lim_{j \to \infty}
                \frac{1}{a}
                h_{F_{r_j}}(0,a)
            =
                \frac{1}{a}
                h_{F_\infty}(0,a)
        \end{align*}
        uniformly on the interval $[s,r]$. By the arbitrariness of $s$ and $r$, we conclude
        \begin{align*}
            \mathcal{N}_{F}(x)
            h_{F_{\infty}}(0,a)
            =
            \mathcal{N}_{F_{\infty}}(0,a)
            h_{F_{\infty}}(0,a)
        \end{align*}
        for every $a > 0$. Consequently, $\mathcal{N}_{F}(x) = \mathcal{N}_{F_{\infty}}(0,a)$ for every $a \in \{ h_{F_{\infty}}(0,\cdot) > 0 \}$. Theorem~\ref{theorem-the-main-frequency-structure-result} yields that
        \begin{align*}
            F_{\infty}(rv) = r^{ \mathcal{N}_{F}(x) +1 } h(v)
        \end{align*}
        and
        \begin{align*}
            \Delta_{ \S^{n-1} } h(v)
            =
            -
            (n-1)\left( \mathcal{N}_{F}(x) +1 \right)^2 h(v)
        \end{align*}
        for $v \in \S^{n-1}$; the normalization for $|h|$ follows from Theorem~\ref{theorem-the-main-frequency-structure-result} using \eqref{equation-blow-up-analysis-normalization}.
    \end{proof}
    In the following subsection, we show that a suitable subclass of homogeneous spherical harmonics can be used to construct radially homogeneous minimal conformal immersions. In a sense, this should be understood as a converse of Theorem~\ref{theorem-the-main-frequency-structure-result}.
    \begin{remark}
        A simple application of Theorem~\ref{theorem-normalized-blow-up} implies that every element of the family
        \begin{align*}
            \mathcal{F} 
            &= \{ F_{\infty}  \colon F_{\infty} = \lim_{j} F_{r_j} \text{ locally uniformly for some $r_j \rightarrow 0^{+}$} \}
        \end{align*}
        satisfies the conclusions of Theorem~\ref{theorem-normalized-blow-up}.
    \end{remark}

    \subsection{From homogeneous spherical harmonics to conformal immersions}\label{section-from-takahashi-to-conformal-immersions}
    Let $h \colon \S^{n-1} \rightarrow \S^{m-1}(0,\xi)$ be a minimal isometric immersion. Then, by Takahashi \cite[Theorem~3]{takahashi-1966-minimal-immersions-of-riemannian-manifolds}, there exists a homogeneous polynomial $P \colon \R^n \to \R^{m}$ of degree $d$ such that $P(rv) = r^{d} h(v)$ for $r > 0$ and $v \in \S^{n-1}$, and $\Delta_{ \S^{n-1} } h = -d(d+n-2) h$. Consider the radial stretching map $\varphi_{k} \colon \R^{n} \to \R^{n}$, where $y \mapsto |y|^{k-1}y$ for $k \in (0,\infty)$, and $F_{k} = P \circ \varphi_{k}$.
    \begin{lemma}\label{lemma-stretching-the-polynomial}
        If $k = (\xi d)^{-1}$, then $F_{k}(rv) = r^{ \frac{1}{\xi} } h(v)$ is a conformal immersion on $\R^n \setminus \{0\}$ and continuous on $\R^n$.
    \end{lemma}
    \begin{proof}
        Let $f \colon \R \times \S^{n-1} \to \R^n \setminus \{0\}$ and $g \colon \R^m \setminus \{0\} \to \R \times \S^{m-1}$ be the conformal diffeomorphisms $f(r,v) = e^{r} v$ and $g(x) = ( \log(|x|), x/|x| )$. Let $Q = g \circ F_k \circ f$. It clearly holds that
        \begin{align*}
            Q( r,v )
            &=
            g( e^{kdr}h(v) )
            =
            g\left( e^{kdr} \xi \frac{ h(v) }{ \xi } \right)
            =
            \left( kdr + \log(\xi), \frac{ h(v) }{ \xi } \right).
        \end{align*}
        It holds, therefore, that the differential of $Q$ has a block representation
        \begin{align*}
            (DQ)(r,v)
            =
            \begin{bmatrix}
                kd & 0
                \\
                0 & \xi^{-1}(Dh)(v)
            \end{bmatrix},
        \end{align*}
        using the identifications $T( \R \times \S^{n-1} ) \simeq \R \times T\S^{n-1}$ and $T( \R \times \S^{m-1} ) \simeq \R \times T\S^{m-1}$. Since $h$ is an isometric immersion, we have
        \begin{align*}
            (DQ)^T (DQ)(r,v)
            =
            \begin{bmatrix}
                (kd)^2 & 0
                \\
                0 & \xi^{-2} I_{ T\S^{n-1} }
            \end{bmatrix}.
        \end{align*}
        Now, if $k = 1/(\xi d)$, the conformality of~$Q$ follows.
    \end{proof}
    \begin{remark}\label{remark-stretching-the-polynomial-rescaling}
        It is convenient to consider the rescaled version of Lemma~\ref{lemma-stretching-the-polynomial}. That is, we consider an isometric immersion $h \colon \S^{n-1}(0,c) \to \S^{m-1}$. Then, up to rescaling, the same argument shows that the map $F \colon \R^n \to \R^{m}$ defined by $rv \mapsto r^{c} h(cv)$ for $r>0$, $v \in \S^{n-1}$, is a conformal immersion in $\R^n \setminus \{0\}$.
    \end{remark}

\section{Proofs of the main results}\label{section-proof-of-the-main-results}
    We now combine the blow-up classification from Theorem~\ref{theorem-normalized-blow-up} with the homogeneity results of Section~\ref{section-radially-homogeneous-conformal-curves} to prove the main theorems.
    \subsection{Proofs of Theorems~\ref{theorem-full-rank-set} and \ref{theorem-discreteness}}\label{section-full-rank-and-fibers}
    \begin{proof}[Proof of Theorem~\ref{theorem-full-rank-set}]
        We argue by contradiction. Then there exists a sequence $( x_j )$ in $S_F$ converging to $x \in S_F$ for which $x_j \neq x_i \neq x$ for each $j \neq i$. We let $r_j = |x_j-x|$ and lose no generality in supposing that $( r_j )$ is decreasing and $\overline{B}(x, 2r_j ) \subset \Omega$. Let $F_{r_j}$ be normalized as in Theorem~\ref{theorem-normalized-blow-up}; we lose no generality that $F_{r_{j}} \rightarrow F_{\infty}$ uniformly and in $C^1$ in each compact subset of $\R^n$. In particular, $v_j = (x_j-x)/r_j \in \S^{n-1}$ satisfy $DF_{r_j}(v_j) = 0$. Again, by passing to a subsequence and relabeling, we may suppose that $v_j \rightarrow v$. By the stated $C^1$-convergence, it follows that $DF_{\infty}(v) = \lim_{j} DF_{\infty}(v_j) = \lim_{j} DF_{r_j}(v_j) = 0$. Since $\partial_{r}F_{\infty}(v) = \left( \mathcal{N}_{F}(x)+1 \right) h_{\infty}(v)$, we conclude $h_{\infty}(v) = 0$. This contradicts Theorem~\ref{theorem-normalized-blow-up}~(1).
    \end{proof}

    \begin{proof}[Proof of Theorem~\ref{theorem-discreteness}]
        We argue by contradiction. Then there exists $w \in \R^m$ and a sequence $( x_j )$ in $F^{-1}(w)$ converging to $x \in \Omega$ for which $x_j \neq x_i \neq x$ for each $j \neq i$. By passing to a limit as in the proof of Theorem~\ref{theorem-full-rank-set}, we obtain a map $F_{\infty}$ satisfying $0 = F_{\infty}(0) = F_{\infty}(v) = h_{\infty}(v)$ for some $v \in \S^{n-1}$. This contradicts Theorem~\ref{theorem-normalized-blow-up}~(1).
    \end{proof}

    \subsection{Proof of Theorem~\ref{theorem-singular-set-is-empty-if-and-only-if-radial-homogeneity}}\label{section-cone-property}
    We first prove the following theorem.
    \begin{theorem}\label{theorem-recognizing-smooth-points}
        For every $x \in S_F$, every blow-up obtained from Theorem~\ref{theorem-normalized-blow-up} is a non-linear radially homogeneous conformal curve. Moreover, every blow-up at $x \in U_F$ obtained from Theorem~\ref{theorem-normalized-blow-up} is non-constant and linear.
    \end{theorem}
    \begin{proof}
        If $x \in U_F$, the smoothness of $F$ at $x$ implies that
        \begin{align*}
            c_{r}(x) = \frac{ |DF|(x) }{ \sqrt{n} } (1+o(1))
            \quad\text{and}\quad
            \mathcal{N}_{F}(x,r) = o(1).
        \end{align*}
        Hence, any blow-up as in the proof of Theorem~\ref{theorem-normalized-blow-up} is linear, and, in fact, equal to the linear map
        \begin{align*}
            y \mapsto \sqrt{n} \frac{ DF(x)[y] }{ |DF|(x) }.
        \end{align*}
        On the other hand, if $x \in S_F$, the radial homogeneity degree of the blow-up is $\mathcal{N}_{F}(x)+1$ by Theorem~\ref{theorem-normalized-blow-up}. Since
        \begin{align*}
            |DF_{\infty}(0)|
            =
            \lim_{j}
            |DF_{r_j}|(0),
        \end{align*}
        we conclude that $\mathcal{N}_{F}(x) > 0$ if and only if $|DF|(x) = 0$, i.e., $x \in S_F$. The claim follows.
    \end{proof}
    \begin{proof}[Proof of Theorem~\ref{theorem-singular-set-is-empty-if-and-only-if-radial-homogeneity}]
        Here (1) and (2) are equivalent by Theorem~\ref{theorem-recognizing-smooth-points}. Indeed, if $F_{\infty}(rv) = r^{\gamma+1} h(v)$ is a conformal curve, then $S_{F_\infty} \neq \emptyset$ if and only if $\gamma > 0$. To finish, '(3) $\Rightarrow$ (2)' is immediate from the characterization in Theorem~\ref{theorem-the-main-frequency-structure-result} by rescaling the radially homogeneous curve. Similarly, '(2) $\Rightarrow$ (3)' follows from Lemma~\ref{lemma-stretching-the-polynomial}.
    \end{proof}
    \begin{remark}
        In the equidimensional case, $U_F = \Omega$ for every non-constant conformal curve $F \colon \Omega \to ( \R^n, \omega )$ for $\omega = \vol_{\R^n}$ or $-\vol_{\R^n}$. Indeed, we verify Theorem~\ref{theorem-singular-set-is-empty-if-and-only-if-radial-homogeneity} (3). To see this, any isometric immersion $h \colon \S^{n-1}(0,c) \rightarrow \S^{n-1}$ is surjective because $h$ is a closed and open map. In fact, $h$ defines a covering map and, as $\S^{n-1}$ is simply connected, a homeomorphism. Therefore, $h$ is an isometry. Necessarily $c = 1$, so the claim follows.
    \end{remark}
    \subsection{Proof of Theorem~\ref{theorem-singular-set-can-be-non-empty}}\label{section-on-a-coassociative-cone}
    Using Lemma~\ref{lemma-stretching-the-polynomial} (and Remark~\ref{remark-stretching-the-polynomial-rescaling}), we prove Theorem~\ref{theorem-singular-set-can-be-non-empty}. For the statement, we use the normalizations of associative $\omega_{\assoc} \in \Lambda^3 \R^7$ and coassociative $\omega_{\coassoc} \in \Lambda^{4} \R^7$ calibrations from \cite[Definition 12.1.1]{joyce-2007-riemannian-holonomy-groups-and-calibrated-geometry}; the calibrations were introduced in \cite{harvey-lawson-1982-calibrated-geometries}. For the convenience of the reader, we give the precise definitions, using the notation $dx^{ijk} = dx^i \wedge dx^j \wedge dx^{k}$ and $dx^{ijkl} = dx_{ijk} \wedge dx^{l}$,
    \begin{align*}
        \omega_{\assoc}
        &=
        dx^{123}
        +
        dx^{145}
        +
        dx^{167}
        +
        dx^{246}
        -
        dx^{257}
        -
        dx^{347}
        -
        dx^{356},
        \\
        \omega_{\coassoc}
        &=
        dx^{4567}
        +
        dx^{2367}
        +
        dx^{2345}
        +
        dx^{1357}
        -
        dx^{1346}
        -
        dx^{1256}
        -
        dx^{1247}.
    \end{align*}
    We first consider the coassociative case. We recall that Ejiri constructs in \cite{ejiri-1981-totally-real-submanifolds-in-a-6-sphere} an example of an isometric minimal immersion $h \colon \R^{4} \supset \S^{3}\left(0,4\right) \to \S^{6} \subset \R^7$ such that the cone over the image of $h$ is calibrated by the coassociative calibration $\omega_{\coassoc}$; indeed, $h$ is totally real and minimal, so the cone is calibrated by $\omega_{coassoc}$ (see e.g.~\cite[page~258]{joyce-2007-riemannian-holonomy-groups-and-calibrated-geometry}). In light of this, Lemma~\ref{lemma-stretching-the-polynomial} (and Remark~\ref{remark-stretching-the-polynomial-rescaling}) implies that the map $F \colon \R^4 \to ( \R^7, \omega_{\coassoc} )$, where $F(rv) = r^{4} h(4v)$, defines a conformal curve. Given that $h$ is real-analytic, the radial scaling $r^4$ implies that $F$ is real-analytic. The fact that $h$ is a $24$-to-$1$ map is established in \cite[Remark~(3), pages~580-581]{dillen-verstraelen-vrancken-1990-classification-of-totally-real-3-dimensional-submanifolds-of-s-6}. In particular, $F$ is $24$-to-$1$ outside the origin.

    We move to the associative case. We describe a minimal isometric immersion $h \colon \R^3 \supset \S^{2}\left( 0, \sqrt{6} \right) \rightarrow \S^{6} \subset \R^7$ known as the Bor{\r{u}}vka--Veronese sphere. By Sekigawa \cite[Example 2]{sekigawa-1983-almost-complex-submanifolds-of-a-6-dimensional-sphere}, the immersion satisfies $$h(\sqrt{6}v) = c ( Y_1(v),\dots,Y_7(v) ), \quad v \in \S^{2},$$ where $Y_1,\dots, Y_7$ form an $L^2$-orthonormal basis of spherical harmonics of degree three. The fact that these describe the Bor{\r{u}}vka--Veronese sphere follows from Calabi \cite{calabi-1967-minimal-immersions-of-surfaces-in-euclidean-spaces}; see also \cite{do-carmo-wallach-1971-minimal-immersions-of-spheres-into-spheres}. The cone over a two-dimensional almost complex submanifold of $\S^{6}$ is known to be associative (see e.g.~\cite[page 187]{bryant-1982-submanifolds-and-special-structures-on-the-octonians} or \cite[page~258]{joyce-2007-riemannian-holonomy-groups-and-calibrated-geometry} for further discussion). Hence $$F \colon \R^3 \to ( \R^7, \omega_{\assoc} ), \quad rv \mapsto r^{ \sqrt{6} } h( \sqrt{6}v )$$ is a conformal curve. The radial homogeneity implies $F \in C^{2}( \R^3, \R^7 ) \setminus C^{3}( \R^3, \R^7 )$. 
    
    It remains to verify the injectivity of $F$, or, equivalently, the injectivity of $h$. It suffices to prove that
        \begin{align*}
            \langle h(\sqrt{6}v), h(\sqrt{6}w) \rangle = 1
            \quad\text{if and only if}\quad
            v = w.
        \end{align*}
        By the addition theorem for spherical harmonics, see e.g. \cite[Theorem~2, page 10]{muller-1966-spherical-harmonics}, the Legendre polynomial $P_{3}$ of degree three satisfies
        \begin{align*}
            1
            =
            c^{2}
            \sum_{j=1}^{7} Y_{j}(v) Y_{j}(w)
            =
            \frac{ 7 }{ 4 \pi }
            c^2
            P_{3}( \langle v, w \rangle ),
        \end{align*}
        where explicitly
        \begin{align*}
            P_{3}( t ) = \frac{1}{2}\left( 5 t^3 - 3 t \right).
        \end{align*}
        First, the normalization of $h(\sqrt{6}v)$ having norm one for every $v \in \S^{2}$ guarantees
        \begin{align*}
            \frac{ 7 }{ 4 \pi }
            c^2
            =
            1.
        \end{align*}
        Second, since $P_{3}(t) = 1$ for $t \in [-1,1]$ if and only if $t = 1$, the claimed injectivity of $F$ follows.
\begin{remark}\label{remark-cayley-is-also-destroyed}
    The Cayley calibration $\omega_{\Cayley} \in \Lambda^4 \R^8$ was introduced by Harvey and Lawson in \cite{harvey-lawson-1982-calibrated-geometries}. For the definition of the calibration, we follow the normalization \cite[equation (12.10)]{joyce-2007-riemannian-holonomy-groups-and-calibrated-geometry}. More precisely, using the decomposition $\R^8 = \R \times \R^{7}$, we consider the coordinate projections $\pi_1 \colon \R^8 \to \R$ and $\pi_2 \colon \R^8 \to \R^7$. Then $$\omega_{\Cayley} = \pi_{1}^{*}\vol_{\R} \wedge \pi_{2}^{*}\omega_{\assoc} + \pi_{2}^{*}\omega_{\coassoc}.$$
    By postcomposing the coassociative curve from Theorem~\ref{theorem-singular-set-can-be-non-empty} with the inclusion $\R^7 \simeq \{0\} \times \R^7 \subset  \R^8$, we obtain a radially homogeneous conformal curve $\R^4 \to ( \R^8, \omega_{\Cayley} )$ with a non-empty singular set.
\end{remark}

   \subsection{Proof of Theorem~\ref{theorem-special-lagrangian-spheres}}\label{section-special-lagrangian}
    Before proceeding with the proof, we need further setup. We first recall that the Special Lagrangian calibrations are defined by
    \begin{align*}
        \omega_{\SL,\theta} = \mathrm{Re}\left( e^{i\theta} dz_1 \wedge \dots \wedge dz_n \right),
    \end{align*}
    where $dz_j$ is the standard line element in the $j$'th factor of $\C^n = \C \times \dots \times \C$. The parameter $\theta \in \R$ is related to the \emph{phase} of the Special Lagrangian. Harvey and Lawson proved that a connected Lagrangian submanifold $M \subset \C^n$ is \emph{minimal} if and only if it is calibrated by some $\omega_{\SL,\theta}$ \cite[Proposition~2.17, page 96]{harvey-lawson-1982-calibrated-geometries}.

    We also recall the following well-known equivalence, see e.g.~\cite[Proposition~1.1]{ohnita-2007-stability-and-rigidity-of-special-lagrangian-cones-over-certain-minimal-legendrian-orbits}, \cite[page 45]{castro-li-urbano-2006-hamiltonian-minimal-lagrangian-submanifolds-in-complex-space-forms}, and \cite[Proposition~2.5]{haskins-2004-special-lagrangian-cones}:
    \begin{enumerate}
        \item a map $h \colon \S^{n-1}(0,c) \rightarrow \S^{2n-1}$ is a minimal Legendrian isometric immersion;
        \item the composition $g = \pi \circ h \colon \S^{n-1}(0,c) \rightarrow \C P^{n-1}$ is a minimal Lagrangian isometric immersion, where $\pi \colon \S^{2n-1} \to \C P^{n-1}$ is the Hopf fibration (where $\C P^{n-1}$ is the complex projective space equipped with the standard metric making $\pi$ a Riemannian submersion);
        \item the conformal radial extension $F \colon \R^n \to \C^n$, defined by $F(rv) = r^{c} h(cv)$ for $(r,v) \in [0,\infty) \times \S^{n-1}$, defines a conformal curve $\R^n \to ( \C^n, \omega_{\SL,\theta} )$ for some $\theta \in \R$.
    \end{enumerate}
    It is more standard to formulate (3) using the map $rv \mapsto G(rv) = r h(v)$ but the minimality of $F$ and $G$ are equivalent. The scaling $r^c$ is chosen to guarantee conformality by Lemma~\ref{lemma-stretching-the-polynomial}. Using this choice, (3) is equivalent to the image of $F$ being calibrated by a Special Lagrangian calibration $\omega_{\SL,\theta} \in \Lambda^n \C^n$.

    Now, we start with the proof of Theorem~\ref{theorem-special-lagrangian-spheres}. We aim to verify Theorem~\ref{theorem-singular-set-is-empty-if-and-only-if-radial-homogeneity} (3). That is, we claim that any isometric immersion $h \colon \S^{n-1}(0,c) \rightarrow \S^{2n-1}$ for which the cone $C$ is calibrated by a Special Lagrangian calibration $\omega_{\SL,\theta}$ is such that $C$ is an $n$-plane. To this end, by the equivalence above, Ejiri's theorem \cite[Theorem, page 1]{ejiri-1982-totally-real-minimal-immersions-of-n-dimensional-real-space-forms-into-n-dimensional-complex-space-forms} gives that $g$ is totally geodesic. Since $h$ is a horizontal lift of $g$ over $\pi$ (see e.g.~\cite[page 45]{castro-li-urbano-2006-hamiltonian-minimal-lagrangian-submanifolds-in-complex-space-forms}), it holds that $h$ is totally geodesic. By comparing the sectional curvatures of $\S^{n-1}(0,c)$ and $\S^{2n-1}$, for $h$ to be totally geodesic, we must have $c = 1$. Hence, the discussion in Remark~\ref{remark-spectrum-of-the-laplacian} implies that $F$ is linear. Therefore, Theorem~\ref{theorem-singular-set-is-empty-if-and-only-if-radial-homogeneity} (3) follows for $\omega_{\SL,\theta}$ and the claimed regularity of Special Lagrangian curves follows.
    \subsection*{Acknowledgements}
    We thank Pekka Pankka and Haohang Zhang for discussions on the topics of the manuscript.

    \subsection*{AI declaration}
    Generative AI tools were used during the development and preparation of this manuscript.
    
    In interactions with the author, GPT-6 Sol derived the identities underlying Proposition~\ref{proposition-variance-of-V} and Lemma~\ref{lemma-key-identities-for-the-energy} under the additional assumption that $S_F = \emptyset$.

    Microsoft Copilot was used primarily for critical proofreading, literature searches, checking standard arguments, and improving the exposition.

    The resulting calculations and arguments were checked, revised, and placed in their present form by the author. The author takes full responsibility for the content of the manuscript.

\providecommand{\bysame}{\leavevmode\hbox to3em{\hrulefill}\thinspace}
\providecommand{\MR}{\relax\ifhmode\unskip\space\fi MR }
\providecommand{\MRhref}[2]{%
  \href{http://www.ams.org/mathscinet-getitem?mr=#1}{#2}
}
\providecommand{\href}[2]{#2}


\begin{thebibliography}{{Smi}11}

\bibitem[Alm79]{almgren-1979-dirichlets-problem-for-multiple-valued-functions-and-the-regularity-of-mass-minimizing-integral-currents}
Frederick~J. Almgren, Jr., \emph{Dirichlet's problem for multiple valued functions and the regularity of mass minimizing integral currents}, Minimal submanifolds and geodesics ({P}roc. {J}apan-{U}nited {S}tates {S}em., {T}okyo, 1977), North-Holland, Amsterdam-New York, 1979, pp.~1--6. \MR{574247}

\bibitem[Alm86]{almgren-1986-optimal-isoperimetric-inequalities}
F.~Almgren, \emph{Optimal isoperimetric inequalities}, Indiana Univ. Math. J. \textbf{35} (1986), no.~3, 451--547. \MR{855173}

\bibitem[BIM25]{broder-iliashenko-madnick-2025-hyperbolicity-and-schwarz-lemmas-in-calibrated-geometry}
Kyle {Broder}, Anton {Iliashenko}, and Jesse {Madnick}, \emph{{Hyperbolicity and Schwarz Lemmas in Calibrated Geometry}}, arXiv e-prints (2025), arXiv:2507.16313.

\bibitem[Bry82]{bryant-1982-submanifolds-and-special-structures-on-the-octonians}
Robert~L. Bryant, \emph{Submanifolds and special structures on the octonians}, J. Differential Geometry \textbf{17} (1982), no.~2, 185--232. \MR{664494}

\bibitem[Cal67]{calabi-1967-minimal-immersions-of-surfaces-in-euclidean-spaces}
Eugenio Calabi, \emph{Minimal immersions of surfaces in {E}uclidean spheres}, J. Differential Geometry \textbf{1} (1967), 111--125. \MR{233294}

\bibitem[CKM20]{cheng-karigiannis-madnick-2020-bubble-tree-convergence-of-conformally-cross-product-preserving-maps}
Da~Rong Cheng, Spiro Karigiannis, and Jesse Madnick, \emph{Bubble tree convergence of conformally cross product preserving maps}, Asian J. Math. \textbf{24} (2020), no.~6, 903--984. \MR{4312754}

\bibitem[CLU06]{castro-li-urbano-2006-hamiltonian-minimal-lagrangian-submanifolds-in-complex-space-forms}
Ildefonso Castro, Haizhong Li, and Francisco Urbano, \emph{Hamiltonian-minimal {L}agrangian submanifolds in complex space forms}, Pacific J. Math. \textbf{227} (2006), no.~1, 43--63. \MR{2247872}

\bibitem[dCW71]{do-carmo-wallach-1971-minimal-immersions-of-spheres-into-spheres}
Manfredo~P. do~Carmo and Nolan~R. Wallach, \emph{Minimal immersions of spheres into spheres}, Ann. of Math. (2) \textbf{93} (1971), 43--62. \MR{278318}

\bibitem[DVV90]{dillen-verstraelen-vrancken-1990-classification-of-totally-real-3-dimensional-submanifolds-of-s-6}
F.~Dillen, L.~Verstraelen, and L.~Vrancken, \emph{Classification of totally real {$3$}-dimensional submanifolds of {$S^6(1)$} with {$K\ge 1/16$}}, J. Math. Soc. Japan \textbf{42} (1990), no.~4, 565--584. \MR{1069845}

\bibitem[Eji81]{ejiri-1981-totally-real-submanifolds-in-a-6-sphere}
Norio Ejiri, \emph{Totally real submanifolds in a {$6$}-sphere}, Proc. Amer. Math. Soc. \textbf{83} (1981), no.~4, 759--763. \MR{630028}

\bibitem[Eji82]{ejiri-1982-totally-real-minimal-immersions-of-n-dimensional-real-space-forms-into-n-dimensional-complex-space-forms}
\bysame, \emph{Totally real minimal immersions of {$n$}-dimensional real space forms into {$n$}-dimensional complex space forms}, Proc. Amer. Math. Soc. \textbf{84} (1982), no.~2, 243--246. \MR{637177}

\bibitem[Geh73]{gehring-1973-the-lp-integrability-of-the-partial-derivatives-of-a-quasiconformal-mapping}
F.~W. Gehring, \emph{The {$L\sp{p}$}-integrability of the partial derivatives of a quasiconformal mapping}, Acta Math. \textbf{130} (1973), 265--277. \MR{402038}

\bibitem[GT01]{gilbarg-trudinger-2001-elliptic-partial-differential-equations-of-second-order}
David Gilbarg and Neil~S. Trudinger, \emph{Elliptic partial differential equations of second order}, Classics in Mathematics, Springer-Verlag, Berlin, 2001, Reprint of the 1998 edition. \MR{1814364}

\bibitem[Has04]{haskins-2004-special-lagrangian-cones}
Mark Haskins, \emph{Special {L}agrangian cones}, Amer. J. Math. \textbf{126} (2004), no.~4, 845--871. \MR{2075484}

\bibitem[HL82]{harvey-lawson-1982-calibrated-geometries}
Reese Harvey and H.~Blaine Lawson, Jr., \emph{Calibrated geometries}, Acta Math. \textbf{148} (1982), 47--157. \MR{666108}

\bibitem[HL87]{hardt-lin-1987-mappings-minimizing-the-lp-norm-of-the-gradient}
Robert Hardt and Fang-Hua Lin, \emph{Mappings minimizing the {$L^p$} norm of the gradient}, Comm. Pure Appl. Math. \textbf{40} (1987), no.~5, 555--588. \MR{896767}

\bibitem[HPP23]{heikkila-pankka-prywes-2023-quasiregular-curves-of-small-distortion-in-product-manifolds}
Susanna Heikkil\"{a}, Pekka Pankka, and Eden Prywes, \emph{Quasiregular curves of small distortion in product manifolds}, J. Geom. Anal. \textbf{33} (2023), no.~1, Paper No. 1, 44. \MR{4502762}

\bibitem[IK25]{iliashenko-karigiannis-2025-a-special-class-of-k-harmonic-maps-inducing-calibrated-fibrations}
Anton Iliashenko and Spiro Karigiannis, \emph{A special class of {$k$}-harmonic maps inducing calibrated fibrations}, Math. Res. Lett. \textbf{32} (2025), no.~5, 1481--1519. \MR{5009217}

\bibitem[IKM26]{iliashenko-karigiannis-madnick-2026-montels-theorem}
Anton {Iliashenko}, Spiro {Karigiannis}, and Jesse {Madnick}, \emph{{Montel's theorem and tautness in calibrated geometry}}, arXiv e-prints (2026), arXiv:2606.31393.

\bibitem[Iko26a]{ikonen-2025-entire-conformal-curves}
Toni Ikonen, \emph{Entire conformal curves are affine or have super-{E}uclidean energy growth}, Proc. Amer. Math. Soc. \textbf{154} (2026), no.~9, 3837--3847. \MR{5096373}

\bibitem[Iko26b]{ikonen-2026-pushforward-of-currents-under-sobolev-maps}
Toni Ikonen, \emph{Pushforward of currents under {S}obolev maps}, Journal of Functional Analysis \textbf{290} (2026), no.~12, 111441.

\bibitem[{Ili}26]{iliashenko-2026-branch-obstructions-for-area-minimizing-currents-and-flat-critical-points-of-weakly-conformal-maps}
Anton {Iliashenko}, \emph{{Branch Obstructions for Area-Minimizing Currents and Flat Critical Points of Weakly-Conformal Maps}}, arXiv e-prints (2026), arXiv:2611.12461.

\bibitem[IP24]{ikonen-pankka-2024-liouvilles-theorem-in-calibrated-geometries}
Toni {Ikonen} and Pekka {Pankka}, \emph{{Liouville's theorem in calibrated geometries}}, arXiv e-prints (J. Anal. Math., to appear) (2024), arXiv:2410.02722.

\bibitem[IP26]{ikonen-pim-2026-a-rescaling-principle-for-quasiregular-curves-with-applications-to-hyperbolicity}
Toni {Ikonen} and Jonathan {Pim}, \emph{{A rescaling principle for quasiregular curves with applications to hyperbolicity}}, arXiv e-prints (2026), arXiv:2607.09353.

\bibitem[Iwa98]{iwaniec-1998-the-gehring-lemma}
Tadeusz Iwaniec, \emph{The {G}ehring lemma}, Quasiconformal mappings and analysis ({A}nn {A}rbor, {MI}, 1995), Springer, New York, 1998, pp.~181--204. \MR{1488451}

\bibitem[Joy07]{joyce-2007-riemannian-holonomy-groups-and-calibrated-geometry}
Dominic~D. Joyce, \emph{Riemannian holonomy groups and calibrated geometry}, Oxford Graduate Texts in Mathematics, vol.~12, Oxford University Press, Oxford, 2007. \MR{2292510}

\bibitem[LF76]{lelong-ferrand-1976-geometrical-interpretations-of-scalar-curvature-and-regularity-of-conformal-homeomorphisms}
Jacqueline Lelong-Ferrand, \emph{Geometrical interpretations of scalar curvature and regularity of conformal homeomorphisms}, Differential geometry and relativity, dath. Phys. Appl. Math., vol. Vol. 3, Reidel, Dordrecht-Boston, Mass., 1976, pp.~91--105. \MR{442846}

\bibitem[Lin06]{lindqvist-2006-notes-on-the-p-laplace-equation}
Peter Lindqvist, \emph{Notes on the {$p$}-{L}aplace equation}, Report. University of Jyv\"askyl\"a{} Department of Mathematics and Statistics, vol. 102, University of Jyv\"askyl\"a, Jyv\"askyl\"a, 2006. \MR{2242021}

\bibitem[M{\"u}l66]{muller-1966-spherical-harmonics}
Claus M{\"u}ller, \emph{Spherical harmonics}, Lecture Notes in Mathematics, vol.~17, Springer-Verlag, Berlin-New York, 1966. \MR{199449}

\bibitem[Ohn07]{ohnita-2007-stability-and-rigidity-of-special-lagrangian-cones-over-certain-minimal-legendrian-orbits}
Yoshihiro Ohnita, \emph{Stability and rigidity of special {L}agrangian cones over certain minimal {L}egendrian orbits}, Osaka J. Math. \textbf{44} (2007), no.~2, 305--334. \MR{2351004}

\bibitem[Pan20]{pankka-2020-quasiregular-curves}
Pekka Pankka, \emph{Quasiregular curves}, Ann. Acad. Sci. Fenn. Math. \textbf{45} (2020), no.~2, 975--990. \MR{4112271}

\bibitem[Sek83]{sekigawa-1983-almost-complex-submanifolds-of-a-6-dimensional-sphere}
Kouei Sekigawa, \emph{Almost complex submanifolds of a {$6$}-dimensional sphere}, Kodai Math. J. \textbf{6} (1983), no.~2, 174--185. \MR{702939}

\bibitem[{Smi}11]{smith-2011-a-theory-of-multiholomorphic-maps}
Aaron~M. {Smith}, \emph{{A theory of multiholomorphic maps}}, arXiv e-prints (2011), arXiv:1112.1471.

\bibitem[Tak66]{takahashi-1966-minimal-immersions-of-riemannian-manifolds}
Tsunero Takahashi, \emph{Minimal immersions of {R}iemannian manifolds}, J. Math. Soc. Japan \textbf{18} (1966), 380--385. \MR{198393}

\bibitem[Tak94]{takeuchi-1994-some-conformal-properties-of-p-harmonic-maps-and-a-regularity-for-sphere-valued-p-harmonic-maps}
Hiroshi Takeuchi, \emph{Some conformal properties of {$p$}-harmonic maps and a regularity for sphere-valued {$p$}-harmonic maps}, J. Math. Soc. Japan \textbf{46} (1994), no.~2, 217--234. \MR{1264939}

\bibitem[Uhl77]{uhlenbeck-1977-regularity-for-a-class-of-non-linear-elliptic-systems}
K.~Uhlenbeck, \emph{Regularity for a class of non-linear elliptic systems}, Acta Math. \textbf{138} (1977), no.~3-4, 219--240. \MR{474389}

\end{thebibliography}
\end{document}